\documentclass [11pt,oneside]{article}

\usepackage{amssymb}
\usepackage{amsfonts}
\usepackage{amsmath}
\usepackage{amsthm}
\usepackage{epsfig}
\usepackage{wrapfig}

\newtheorem{lemma}{Lemma}[section]
\newtheorem{teo}[lemma]{Theorem}
\newtheorem{rem}[lemma]{Remark}
\newtheorem{prop}[lemma]{Proposition}
\newtheorem{cor}[lemma]{Corollary}
\newtheorem{conj}[lemma]{Conjecture}

\newcommand{\matN} {\ensuremath {\mathbb{N}}}

\newcommand{\matZ} {\ensuremath {\mathbb{Z}}}

\newcommand{\matH} {\ensuremath {\mathbb{H}}}

\newcommand{\calK} {\ensuremath {\mathcal{K}}}

\newcommand{\calS} {\ensuremath {\mathcal{S}}}

\newcommand{\calD} {\ensuremath {\mathcal{D}}}

\newcommand{\calT} {\ensuremath {\mathcal{T}}}

\newcommand{\calH} {\ensuremath {\mathcal{H}}}
\newcommand{\calO} {\ensuremath {\mathcal{O}}}

\newcommand{\nota} [1] {\caption{\footnotesize{#1}}}

\newcommand{\compo}{\,{\scriptstyle\circ}\,}

\newfont{\Got}{eufm10 scaled 1200}
\newcommand{\permu}{{\hbox{\Got S}}}

\font\titsc=cmcsc10 scaled 1200

\newcommand{\dimo}[1]{\vspace{2pt}\noindent\textit{Proof of \ref{#1}}.\ }
\newcommand{\finedimo}{{\hfill\hbox{$\square$}\vspace{2pt}}}
\newcommand{\mettifig}[1]{\epsfig{file=#1}}
\newcommand{\Ker}{\ensuremath{\mathrm{Ker\,}}}

\newcommand{\hatY}{\widehat{Y}}
\newcommand{\Aut}{\textrm{Aut}}
\newcommand{\bysame}{\leavevmode\hbox to3em{\hrulefill}\thinspace}

\author{Fran\c{c}ois \titsc{Costantino}
\and Roberto \titsc{Frigerio}
\and Bruno \titsc{Martelli}
\and Carlo \titsc{Petronio}}

\title{Triangulations of 3-manifolds, hyperbolic\\
relative handlebodies, and Dehn filling}

\begin{document}

\maketitle

\begin{abstract}
    \noindent
    We establish a bijective correspondence between the set $\calT_n$ of
    3-di\-men\-sio\-nal triangulations with $n$ tetrahedra
    and a certain class $\calH_n$ of
    relative handlebodies (\emph{i.e.}~handlebodies with boundary loops,
    as defined by Johannson)
    of genus $n+1$.

    We show that the manifolds in $\calH_n$ are hyperbolic (with geodesic boundary, and cusps corresponding to the loops),
    have least possible volume, and simplest boundary loops.

    Mirroring the elements of $\calH_n$ in their geodesic boundary we obtain
    a set $\calD_n$ of cusped hyperbolic manifolds, previously considered by D.~Thurston
    and the first named author. We show that also $\calD_n$ corresponds bijectively to $\calT_n$, and we
    study some Dehn fillings of the manifolds in $\calD_n$.
    As consequences of our constructions, we also show that:
    \begin{itemize}
    \item A triangulation of
    a 3-manifold is uniquely determined up to isotopy by its 1-skeleton;
    \item If a 3-manifold $M$ has an ideal triangulation with edges of valence at least $6$, then $M$
      is hyperbolic and the edges are homotopically non-trivial, whence homotopic to geodesics;
    \item Every finite group $G$ is the isometry group of a closed hyperbolic 3-manifold with
      volume less than ${\rm const}\times |G|^9$.
    \end{itemize}
  \vspace{4pt}

\end{abstract}

\section*{Introduction}
Let us extend the traditional notion of triangulation
of a 3-manifold by calling \emph{triangulation} any combinatorial pattern of
face-pairings between a finite number of tetrahedra.
So a triangulation may actually not define a manifold. We also define
a \emph{hyperbolic} manifold to be a complete finite-volume hyperbolic 3-manifold
with (possibly empty) geodesic boundary.

In this paper we analyze the class $\calT$ of triangulations
in geometric terms, showing that it can be identified to
a set $\calH$ of hyperbolic 3-manifolds, and to the set $\calD$ of their ``doubles,''
already considered in~\cite{Co:Th} by D.~Thurston and the first named author.
The correspondences we construct shade new light on each of the objects
we consider and allow us to prove several topological and geometric facts.

We now state our main results. 
Proofs will be given in Sections~\ref{tria:section} to~\ref{isom:section}.

\paragraph{Relative handlebodies}
Following Johannson~\cite{johannson}, a \emph{relative handlebody}
$(H,\Gamma)$ is a (possibly non-orientable)
handlebody $H$ with a finite system $\Gamma$ of disjoint
loops on $\partial H$.
A \emph{hyperbolic structure} on $(H,\Gamma)$ is
a finite-volume complete hyperbolic structure with totally geodesic
boundary on $H\setminus\Gamma$, so each loop in $\Gamma$ geometrically
corresponds to a cusp based on a strip (an annulus or a M\"obius band).

We call \emph{pant-meridinal complexity} of $\Gamma$ the minimum $c(\Gamma)$ of
$|\Gamma\cap\partial D|$ where $D$ is a system of disjoint properly
embedded discs in $H$ such that $\partial D$ cuts $\partial H$ into a
union of pairs of pants.

\begin{figure}
\begin{minipage}{.03\textwidth}\hfil
\end{minipage}
\begin{minipage}{.3\textwidth}
\centering
\mettifig{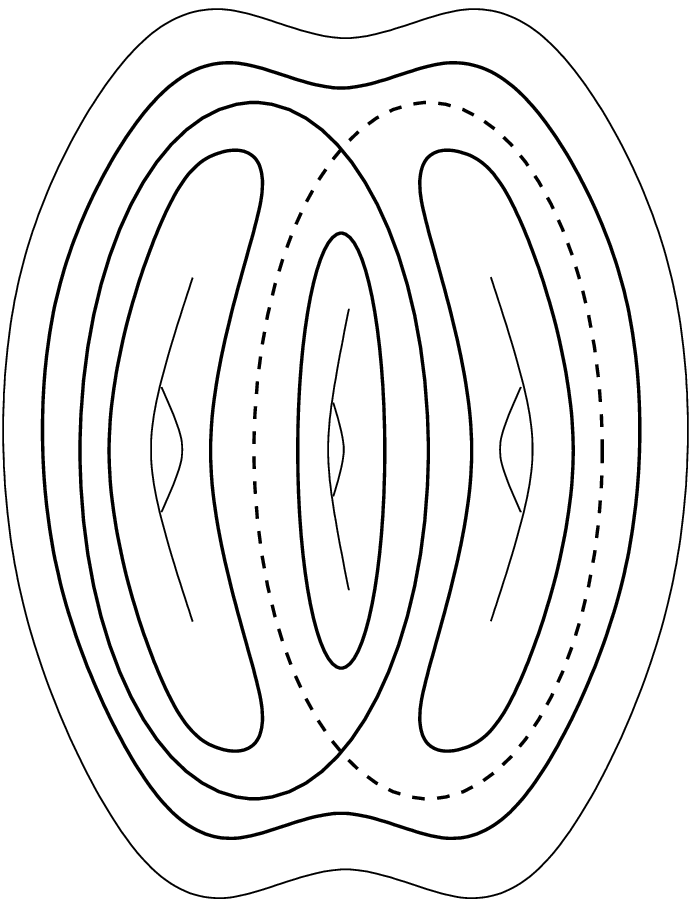, width = .9\textwidth}
\vglue -4mm
\end{minipage}
\begin{minipage}{.03\textwidth}\hfil
\end{minipage}
\begin{minipage}{.57\textwidth}
{If $T$ is the triangulation of $S^3$ obtained by gluing two tetrahedra
along the identity of their boundary, the corresponding genus-$3$
relative handlebody $N(T)=(H,\Gamma)$ is as shown.
Its geometric realization, which has $6$ annular cusps and decomposes into two
regular ideal octahedra, is also the building block
of the elements of $\calD$ (defined below).
}
\end{minipage}
\nota{An example of $N(T)$.}
\label{example:fig}
\end{figure}

\paragraph{From triangulations to relative handlebodies}
If $T\in\calT$ we now define a relative handlebody $N(T)$ as follows.
First, we glue together the tetrahedra of $T$ according to the face-pairings,
getting the support $|T|$ of $T$.
Then, we remove from $|T|$ an open neighbourhood of the vertices, getting a
space $M(T)$. Next, we remove from $M(T)$ an
open neighbourhood $W$ of the edges, which gives a handlebody $H$.
Last, we note that $\partial W$ consists of strips,
we define $\Gamma$ as the system of cores of $\partial W$, and we set
$N(T)=(H,\Gamma)$. We often denote by $N(T)$
also the non-compact manifold $H\setminus\Gamma$, the context
making clear which definition of $N(T)$ we are referring to.
An example of $N(T)$ is shown in Fig.~\ref{example:fig}.

In Section~\ref{tria:section} we give $N(T)$ a hyperbolic structure,
by taking a regular hyperbolic ideal  octahedron for each tetrahedron of $T$,
and gluing these octahedra appropriately.
We then establish the following result.
We denote by $\calT_n$ the set of triangulations involving $n$ tetrahedra.

\begin{teo}\label{T:to:H:intro:teo}
The map $T\mapsto N(T)$ just defined gives a bijection
between $\calT_n$ and a set $\calH_n$ of hyperbolic relative handlebodies of genus $n+1$.
Among the hyperbolic $(H,\Gamma)$'s of genus $n+1$, the set $\calH_n$ consists of:
\begin{itemize}
\item those having minimal complexity $c(\Gamma)$, equal to $10\cdot n$;
\item those having minimal volume, equal to $n\cdot v_O$, where
$v_O\approx 3.66386$ is the volume of a hyperbolic regular ideal octahedron.
\end{itemize}
\end{teo}

The map $N:\calT\to\calH:=\bigcup_{n=1}^{\infty} \calH_n$
is a bijection because the decomposition of $N(T)$ into octahedra is
Kojima's canonical one~\cite{Kojima:deco}.
A consequence of this bijection is the following result, which appears
to be new and to have independent interest.
If $Y$ is a compact 3-manifold, we call $T$ a
\emph{triangulation} of $Y$ if a space homeomorphic to
$Y$ can be obtained from $|T|$ by removing open
stars of some of the vertices. Therefore we have an ideal
vertex for each component (if any) of $\partial Y$, and possibly some internal vertices.
When we remove stars of all the vertices, we say that $T$ is an
\emph{ideal} triangulation of $Y$, and $Y$ coincides with the space $M(T)$
already defined above.

\begin{teo}\label{triangulations:determined:intro:teo}
Let $T_0$ and $T_1$ be triangulations of the same compact $3$-manifold $Y$.
Assume that the  $1$-skeleta of $T_0$ and $T_1$ are the same graph contained in $Y$.
Then $T_0$ and $T_1$ are isotopic relatively to the 1-skeleton.
\end{teo}

\paragraph{Tangles}
We now call \emph{tangle} in a compact 3-manifold with boundary
a finite union of disjoint properly embedded arcs.
As above, the complement of an open tubular neighbourhood $W$ of a tangle
is a manifold with boundary loops given by the cores of $\partial W$.
Since every tangle is contained in the 1-skeleton of an ideal triangulation~\cite{gennaro},
we have the following result.

\begin{cor}\label{tangles:enlarged:intro:cor}
Every tangle in every compact $3$-manifold is contained in a tangle whose complement lies in $\calH$
(in particular, it is hyperbolic).
\end{cor}

\paragraph{Mirroring}
If $(H,\Gamma)$ is a relative handlebody we now define its ``double'' $D(H,\Gamma)$ as
follows. We take the orientation covering $\widetilde H$ of $H$ with deck involution $\tau$ (so
$H=\widetilde H/_\tau$), we consider the pre-image $\widetilde\Gamma$ of $\Gamma$ in $\widetilde H$,
and we define $D(H,\Gamma)=\widetilde H/_{\tau_\partial}$, where $\tau_\partial$ is the
restriction of $\tau$ to $\partial\widetilde H\setminus\widetilde\Gamma$.
When $H$ is orientable, $D(H,\Gamma)$ is just $H\setminus\Gamma$ mirrored along its boundary
$(\partial H)\setminus\Gamma$. In general, $D(H,\Gamma)$
is an open orientable manifold with cusps based on tori. If $T\in\calT_n$ we can construct
$D(T):=D(N(T))$, and
define $\calD_n$ as the set of manifolds of this type. We also set
$\calD=\bigcup_{n=1}^\infty \calD_n$.

The manifold $D(T)$ is constructed in~\cite{Co:Th}
as the boundary of the 4-dimensional thickening of a certain 2-dimensional 
polyhedron $P_0(T)$ ``dual'' to $T$
(this construction is not needed here and is summarized in the Appendix).
The analogue of Corollary~\ref{tangles:enlarged:intro:cor} for $\calD$ is proved in~\cite{Co:Th}, and it
says that every closed orientable
3-manifold is a Dehn filling of a manifold in $\calD$,
\emph{i.e.}~\emph{the class $\calD$ is universal for $3$-manifolds under Dehn filling}.
The following result shows that the elements of
$\calD$ correspond bijectively to those of $\calT$, so they can be easily classified
in purely combinatorial terms.

\begin{teo}\label{T:to:D:intro:teo}
Every member of $\calD$ is an orientable cusped
hyperbolic manifold without boundary.
The correspondence $T\mapsto D(T)$ defines a bijection
between $\calT$ and $\calD$.
\end{teo}

We mention here that the set $\calD$ consists of the manifolds which
decompose into ``blocks''
homeomorphic to that shown in Fig.~\ref{example:fig}.
This block is also used
in~\cite{Agol:K-blocks} to construct a hyperbolic manifold from a path in the curve complex of a surface.
In Section~\ref{surg:section}
we show that every manifold in $\calD$ (with a very few exceptions)
contains a closed incompressible genus-2 surface.

\paragraph{Isometries}
Since a Dehn filling of a hyperbolic manifold ``typically'' is hyperbolic,
the class $\calD$ provides a powerful method to construct closed hyperbolic
manifolds. As an application of this method we establish in Section~\ref{isom:section}
the following result, the first statement of which was already known~\cite{Kojima:isom}.

\begin{teo}\label{all:groups:intro:teo}
If $G$ is a finite group then there exists a closed orientable hyperbolic
$3$-manifold $V$ such that the isometry group of $V$ is isomorphic to $G$.
Moreover $V$ can be chosen so that ${\rm vol}(V)\leqslant c\cdot |G|^9$, where
$c>0$ is a constant.
\end{teo}

\paragraph{Exceptional slopes}
In Section~\ref{surg:section} we study some Dehn fillings of the manifolds in $\calD$.
To state our main results we recall that, according to Thurston's hyperbolic Dehn filling theorem, on each cusp of a
finite-volume hyperbolic 3-manifold there is only a finite number of slopes filling along
which one gets a non-hyperbolic 3-manifold. These slopes are called \emph{exceptional}, and
a considerable effort has been devoted to understanding them~\cite{gordon}.
If $T$ is a triangulation, the hyperbolic manifold $D(T)$ has a preferred
horospherical cusp section, and each component of this section corresponds
to an edge of $T$. Moreover the valence of the edge gives a lower
bound for the length of the second shortest geodesic on the component.
This fact
and the Agol-Lackenby 6-theorem~\cite{Agol:length-6,lackenby} imply the following:

\begin{prop}\label{exceptional:intro:prop}
If every edge of $T$ has valence at least $7$ then there is at most one exceptional slope on
each cusp of $D(T)$.
\end{prop}

\paragraph{Hyperbolicity from combinatorics}
Suppose a 3-manifold
$M$ has an \emph{ideal} triangulation $T$ such that each edge has valence at least $6$.
An easy
argument shows that $\chi(T)\leqslant 0$, and that $\chi(T)=0$ precisely
when all valences are $6$. Moreover, in the last case,
the boundary of $M$ is a disjoint union of tori and Klein bottles, and
Thurston's hyperbolicity equations for cusped manifolds
have a very simple solution, given by regular ideal tetrahedra.
Analogously, if all edges have one and the same valence $v\geqslant 7$, the hyperbolicity
equations for the geodesic boundary case described in~\cite{FriPe}
have a simple solution, given by regular truncated tetrahedra
with dihedral angles $2\pi/v$.
An argument based on the Agol-Lackenby machinery~\cite{Agol:length-6, lackenby}
allows us to generalize these facts as follows (see Section~\ref{surg:section}):

\begin{prop}
If $M$ has an ideal triangulation $T$ whose edges have valence at least $6$, then $M$ is hyperbolic,
and the edges of $T$ are homotopically non-trivial relative to $\partial M$.
\end{prop}

If all valences are at least $6$, every edge of $T$ is then homotopic to a geodesic.
This shows that the tetrahedra of $T$ themselves are homotopic to straight truncated
ones, so it seems quite natural to ask whether such tetrahedra
are in fact \emph{isotopic} to straight ones, \emph{i.e.} whether $T$ gives a geometric
decomposition of $M(T)$.  This is true when the valences of the edges are all equal to each other,
as discussed above. And we have verified experimentally the same fact
in~\cite{FriMaPe:exp} for all
triangulations involving at most $4$ tetrahedra.
We therefore propose the following conjecture:

\begin{conj}
If $M$ has an ideal triangulation $T$ whose edges have valence at least $6$, then $T$
is realized by hyperbolic partially truncated tetrahedra.
\end{conj}

\paragraph{Acknowledgments} We thank Simon King for communicating to us a
proof, based on traditional cut-and-paste techniques, of
Theorem~\ref{triangulations:determined:intro:teo} for the special case
of genuine triangulations (without multiple and self-adjacencies).
F.~C.~thanks Dylan Thurston for his
continuous encouragement and illuminating observations.
B.~M.~thanks the Mathematics Department of the University of Austin
for hospitality during May 2003.
C.~P.~thanks the Universit\'e Paris 7 and the Institut de
Math\'ematiques de Jussieu for hospitality during May 2003.

We thank the referee for his very useful suggestions.

\section{Triangulations and\\ hyperbolic relative handlebodies}\label{tria:section}
In this section we give a precise definition of
what we call a triangulation, and we show that the set of
triangulations corresponds to a class of hyperbolic structures on
relative handlebodies. We also provide geometric and
topological information on these handlebodies, including two
intrinsic definitions of their class.

\paragraph{Triangulations}
We call \emph{triangulation} a pair $(\{\Delta_i\}_{i=1}^n,\{g_j\}_{j=1}^{2n})$,
where $n\in\matN$ is positive, each $\Delta_i$ is a copy of the standard
tetrahedron, and the $g_j$'s give a complete system of
simplicial pairings between the faces of the $\Delta_i$'s, such that the gluing of
the $\Delta_i$'s along the $g_j$'s is connected.
Note that this gluing may actually not be a 3-manifold, because the link of a vertex
could be any surface, and the link of the midpoint of an edge could be
a projective plane. Even when the gluing is a manifold, the
$\Delta_i$'s give a triangulation only in a loose sense, because
multiple adjacencies and self-adjacencies are allowed.
Triangulations are viewed up to combinatorial equivalence, and the set
of equivalence classes is denoted by $\calT$.

If $\Delta$ is the tetrahedron, we denote now by $\Delta^*$ the polyhedron obtained
by truncating $\Delta$ at the vertices (\emph{i.e.}, formally, by removing
open stars of the vertices). Note that $\partial\Delta^*$ consists
of four truncation triangles and four
``lateral'' hexagons.
If $T=(\{\Delta_i\}_{i=1}^n,\{g_j\}_{j=1}^{2n})$ is a triangulation then
the $g_j$'s give gluing rules for the lateral hexagons of the
$\Delta_i^*$'s, and we denote by $M(T)$ the result, which may or not
be a manifold. In case $M(T)$ is a manifold, $T$ is called an
\emph{ideal triangulation} of $M(T)$. Note that in this case
$\partial M(T)$ is triangulated (in a loose sense) by the truncation triangles of
the $\Delta_i^*$'s.

\paragraph{From triangulations to relative handlebodies}
Fix a triangulation $T=(\{\Delta_i\}_{i=1}^n,\{g_j\}_{j=1}^{2n})$.
If we remove from each $\Delta_i$ an open neighbourhood of the
whole 1-skeleton, we can still glue what is left along the
restrictions of the $g_j$'s. This corresponds to adding
$1$-handles to a disjoint union of $0$-handles, so the result is a
(possibly non-orientable) handlebody $H$.
Of course $H$ can also be obtained from
$M(T)$ by removing an open neighbourhood of the image of the set of
edges, as discussed in the Introduction.
Each component of this neighbourhood is a cyclic gluing of
wedges, which implies that $\partial H\setminus \partial M(T)$ is
a union of strips (annuli and/or M\"obius bands). We can then
consider the system $\Gamma$ of the cores of these strips, and get
a relative handlebody $(H,\Gamma)=:N(T)$. Note that $M(T)$ is a
manifold if and only if a regular neighbourhood of $\Gamma$ in
$\partial H$ consists of annuli only (no M\"obius strips).

\paragraph{Hyperbolic structure}We have already defined in the Introduction what
we mean by \emph{hyperbolic structure} on a relative handlebody. We
now explain how to construct one such structure on each $N(T)$.

Given a manifold with a triangulation $T=(\{\Delta_i\}_{i=1}^n,\{g_j\}_{j=1}^{2n})$,
there is a general strategy to algorithmically construct a hyperbolic structure on $N(T)$ using $T$.
This strategy, which dates back to Thurston~\cite{thurston:notes}
and was explained in detail in~\cite{FriPe} for the geodesic boundary case,
amounts to choosing the dihedral angles of each $\Delta_i$ along its edges,
which gives $\Delta_i$ the shape of a hyperbolic polyhedron in $\matH^3$
with some points at infinity, and then requiring that the hyperbolic polyhedra
glue up coherently to give a complete structure.

We do not need to reproduce the details of this strategy here. We
only note that, when trying to construct a hyperbolic structure on
a relative handlebody $N(T)$, all the dihedral angles in $T$ are
actually forced to be $0$. So in this case there is no choice to
make: we only need to analyze what is a tetrahedron with $0$
dihedral angles, and to show that gluing such objects we get a
coherent and complete structure.

\begin{figure}
    \begin{center}
    \mettifig{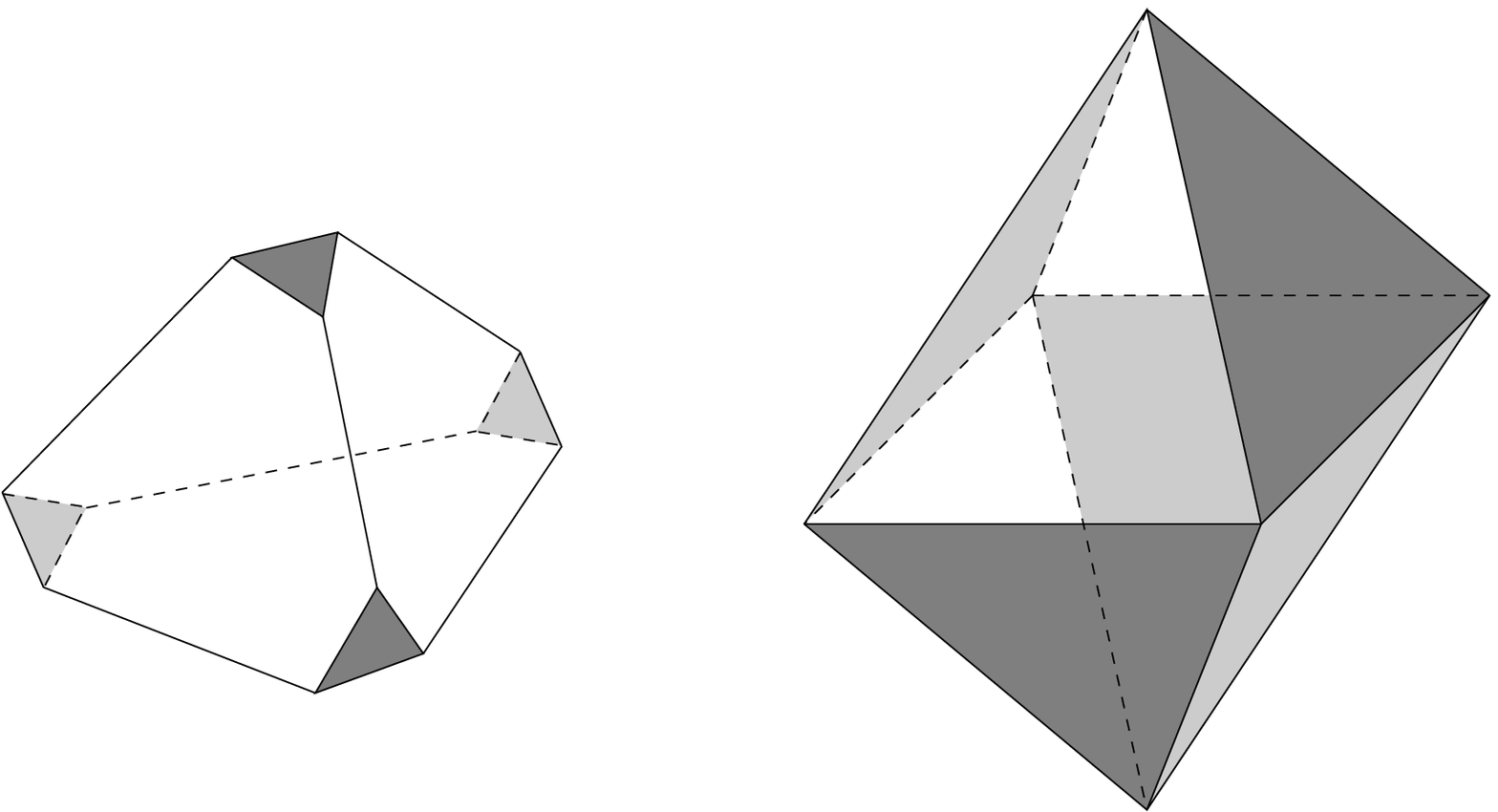,width=8cm}
    \nota{A truncated tetrahedron with dihedral angles zero is
    a regular ideal octahedron.} \label{tetra:octa:fig}
    \end{center}
    \end{figure}

To accomplish the first task, start with a truncated tetrahedron as in
Fig.~\ref{tetra:octa:fig}-left, with black truncation triangles and white
lateral hexagons. Now recall that declaring the dihedral angle
along an edge to be  $0$ geometrically means that the edge disappears into an ideal
vertex. The geometric version of the tetrahedron is therefore an ideal octahedron,
with a checkerboard coloring of the faces, as in Fig.~\ref{tetra:octa:fig}-right.
Each dihedral angle of the octahedron is $\pi/2$, because it is the angle between
a (now ideal) truncation triangle and a lateral ``hexagon,'' now degenerated
into an ideal triangle. It follows that the octahedron is regular. We
denote now by $O$ the regular ideal octahedron with checkerboard coloring
of the faces as in Fig.~\ref{tetra:octa:fig}.

To understand the gluing geometrically, we note that only the
white triangles get glued (in pairs), while the black triangles
give the geodesic boundary. Coherence of the structure along the
gluings is not an issue, because any two ideal triangles are
isometric to each other. Moreover, in the language
of~\cite{FriPe}, there are no internal edges, so the gluing
automatically gives a (possibly incomplete) hyperbolic structure
on $N(T)$. However, completeness is automatic when there are no
ends based on tori~\cite{FriPe}, so indeed $N(T)$ is hyperbolic.
This already implies Corollary~\ref{tangles:enlarged:intro:cor}.

To state our next (easy) result we must recall that if $V$ is
hyperbolic with geodesic boundary and we consider a cusp of $V$
based on a strip, then the strip has a well-defined Euclidean
structure up to rescaling.  If $(H,\Gamma)=N(T)$ and
$\gamma\in\Gamma$ we can then speak of the Euclidean structure of
the strip $U(\gamma)$
on which the cusp corresponding to $\gamma$
is based. We also recall that $\gamma$ comes from an edge of $T$.

\begin{lemma}\label{strip:shape:lem}
If $\gamma$ comes from an edge of valence $q$ then the Euclidean structure
of $U(\gamma)$ is obtained from a rectangle $[0,q]\times[0,1]$ by gluing
$\{0\}\times[0,1]$ to $\{q\}\times[0,1]$.
\end{lemma}

\begin{proof}
The cross-section of the octahedron $O$ at a vertex is a square,
with two opposite black edges and two white edges.  As we construct $U(\gamma)$
we glue together $q$ such squares along white edges.
\end{proof}

\paragraph{Canonical decomposition}
In the projective model of $\matH^3$ consider a polyhedron $P$ with some ideal
and some ultra-ideal vertices, and assume that all the edges of $P$ meet the closure of
$\matH^3$.  Dual to the ultra-ideal vertices of $P$ there are hyperplanes of
$\matH^3$, truncating $P$ along which we get the so-called ``partially truncated polyhedron''
associated to $P$. Kojima~\cite{Kojima:deco} has shown that a finite-volume hyperbolic 3-manifold
$V$ with non-empty geodesic boundary has a canonical decomposition $\calK(V)$ into
partially truncated polyhedra. Combinatorially, $\calK(V)$
corresponds to a gluing $(\{P_i\},\{g_j\})$ of genuine polyhedra ---an
obvious extension of our notion of triangulation, where we allow arbitrary
polyhedra instead of tetrahedra--- and we use the same symbol
$\calK(V)$ to denote both the geometric and the combinatorial version
of the decomposition. When the cusps of $V$ are based on strips only,
$\calK(V)$ is actually dual to the cut-locus of $\partial V$.
We can now establish the following:

\begin{prop}\label{T:is:candeco:prop}
If $T\in\calT$ then $\calK(N(T))$ consists of the regular ideal
octahedra employed to construct the hyperbolic structure. In
particular, $\calK(N(T))$ equals $T$ combinatorially.
\end{prop}

\begin{proof}
Since $\calK(N(T))$ is dual to the cut-locus of $\partial N(T)$,
it suffices to prove that the latter is obtained by gluing
together the cut-loci of the black faces within the individual
octahedra $O$ employed to construct $N(T)$. The cut locus of the
black faces of $O$ is shown in Fig.~\ref{octa:cut:fig}-left
combinatorially, and in Fig.~\ref{octa:cut:fig}-centre
geometrically.
    \begin{figure}
    \begin{center}
    \mettifig{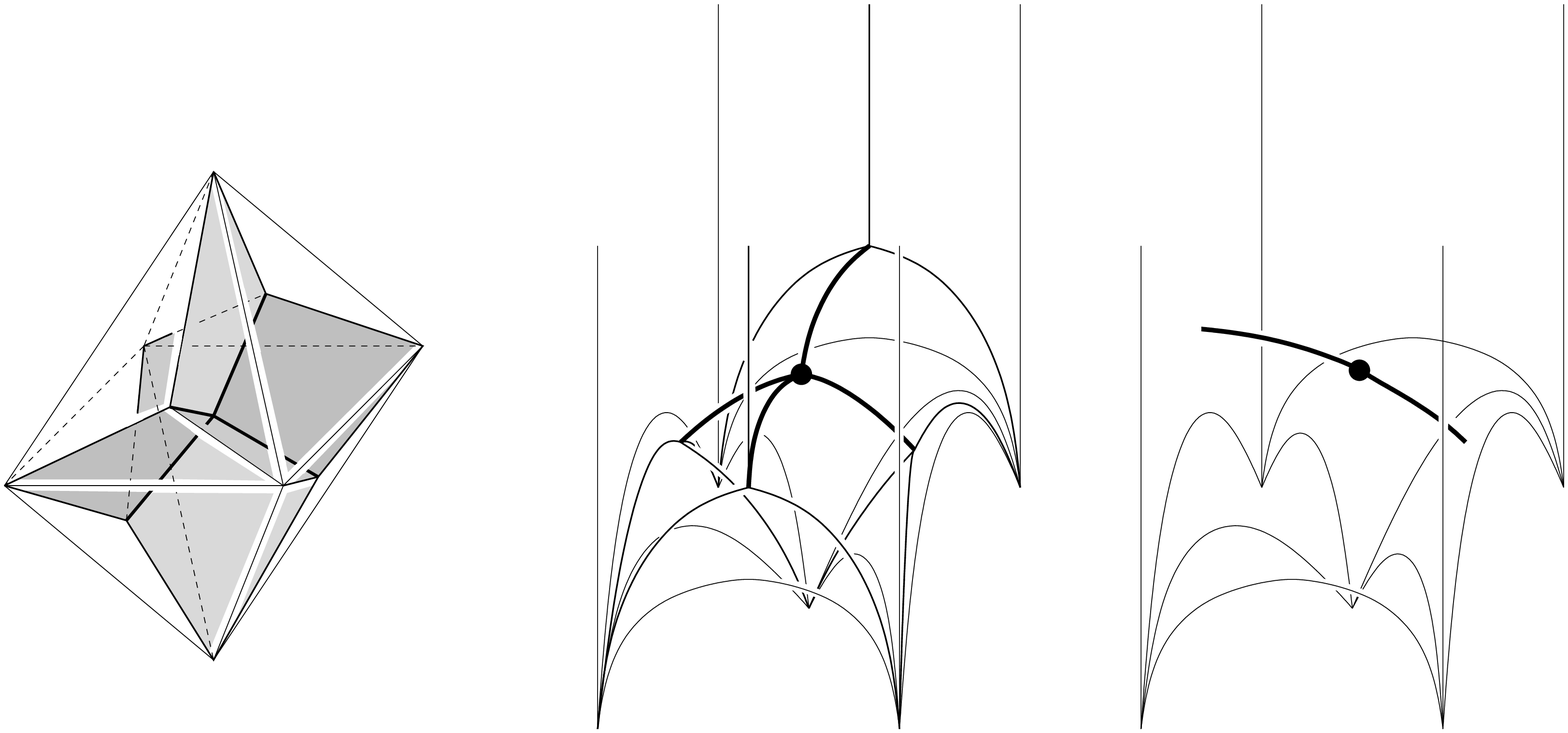,width=12cm}
    \nota{Cut locus of the black faces of $O$, and shortest
    path from a black to a white face.} \label{octa:cut:fig}
    \end{center}
    \end{figure}

To prove our result we must show that if $p\in N(T)$ and $q$ is a
point of $\partial N(T)$ closest to $p$ then $p$ and $q$ belong to
the same copy of $O$. If this is not the case, the shortest
geodesic from $q$ to $p$ starts with an arc which leaves a black
face of a copy of $O$ orthogonally at $q$ and exits from a white
face. The shortest such arc joins the centre of a black face to
the centre of the opposite white face, as shown in
Fig.~\ref{octa:cut:fig}-right. Comparing with
Fig.~\ref{octa:cut:fig}-centre we see that the arc is at least
twice as long as the distance from $p$ to the black faces of the
copy of $O$ in which $p$ lies. This of course gives a
contradiction.
\end{proof}

An alternative proof of the previous proposition could also be given by computing the
\emph{tilts} of the faces of $T$, as explained in~\cite{FriPe}.
The proposition readily implies that the map $T\mapsto
N(T)$ is a bijection between $\calT$ and a set $\calH$ of
hyperbolic relative handlebodies (the main assertion of
Theorem~\ref{T:to:H:intro:teo}). Moreover $\calH$ consists
precisely of the hyperbolic manifolds which decompose into copies
of the octahedron $O$, with black faces on the boundary and white
faces in the interior.

The same result also allows us to prove that a triangulation of a
3-manifold $Y$ is determined by its $1$-skeleton as a subset of
$Y$ (Theorem~\ref{triangulations:determined:intro:teo}). To
state this precisely, denote by $\hatY$ the space obtained from
$Y$ by collapsing each component of $\partial Y$ to a point. Then
note that if we have a triangulation of $Y$ with ideal vertices
at the components of $\partial Y$ and possibly some non-ideal
vertices, then the $1$-skeleton of the triangulation is a graph
contained in $\hatY$, well-defined up to isotopy.

\begin{teo}\label{triangulations:determined:teo}
Let the triangulations $T_0$ and $T_1$ have the same $1$-skeleton $J\subset\hatY$. Then
$T_0$ and $T_1$ are isotopic relatively to $J$.
\end{teo}

\begin{proof}
For $i=0,1$ let $N(T_i)=(H_i,\Gamma_i)$, and
observe that $\hatY\setminus
J\cong H_i\setminus\Gamma_i$, so the identity map
of $\hatY\setminus J$ induces a well-defined homeomorphism
$h:(H_0,\Gamma_0)\to (H_1,\Gamma_1)$.
Moreover, the geometric realization of the tetrahedra
of $T_i$ defines on $H_i\setminus\Gamma_i$ a hyperbolic
metric. Mostow's rigidity theorem then gives a homotopy $\big(h_t:H_0\to H_1\big)_{t\in[0,1]}$
such that $h_0=h$, $h_1$ induces an isometry $H_0\setminus\Gamma_0\to H_1\setminus\Gamma_1$,
and for all $t\in[0,1]$ the following happens:
$$h_t(\partial H_0)=\partial H_1,\quad
h_t(\Gamma_0)=\Gamma_1,\quad h_t^{-1}(\Gamma_1)=\Gamma_0.$$
(These facts are not obvious \emph{a priori}. They follow from the construction
of the homotopy as a convex combination of two maps defined
on the universal cover).
Since $T_i$ is the canonical Kojima decomposition of $H_i\setminus\Gamma_i$, we have $h_1(T_0)=T_1$.
Suppose for a moment that each $h_t$ is a
homeomorphism and $h_t(\Gamma_0)=\Gamma_1$. Then we can use the natural identification $\hatY\setminus
J\cong H_i\setminus\Gamma_i$ and extend $h_t$ to a homeomorphism of $\hatY$ such that
$h_t(J)=J$, thus getting the desired isotopy between $T_0$ and $T_1$.
The rest of the proof is devoted to proving that $\big(h_t\big)_{t\in[0,1]}$
can indeed be replaced by an isotopy mapping $\Gamma_0$ to $\Gamma_1$ for all $t$.

A result of Waldhausen~\cite[Theorem 7.1]{Wald}
implies that there exists an isotopy $\big(h'_t:H_0\to
H_1\big)_{t\in[0,1]}$ such that $h'_0=h_0=h$ and $h'_1=h_1$.
Note that the trace $f'_t$ of $h'_t$ on $\partial H_0$ need not map $\Gamma_0$ to
$\Gamma_1$, so $\big(h'_t\big)_{t\in[0,1]}$ is not the desired isotopy yet.
To modify it, we return to the initial $\big(h_t\big)_{t\in[0,1]}$ and
for $j=0,1$ we define $f_j:\partial H_0\to\partial H_1$ as the restriction of $h_j$.
Now recall that surface homeomorphisms which are homotopic relatively
to the boundary are in fact isotopic~\cite{epstein:acta}. Using
this fact one easily sees that $f_0$ and $f_1$ embed into an isotopy
$\big(f_t:\partial H_0\to\partial H_1\big)_{t\in [0,1]}$ such that
$f_t(\Gamma_0)=\Gamma_1$ for all $t\in [0,1]$.

We have constructed two paths $\big(f'_t\big)_{t\in[0,1]}$ and
$\big(f_t\big)_{t\in[0,1]}$ of homeomorphisms
 $\partial H_0\to\partial H_1$, with the
same ends. Now recall that the space of
isotopically trivial automorphisms
of a closed hyperbolic surface is
simply connected~\cite{earle}. Therefore, if we
denote by $H'_i$ the handlebody obtained from $H_i$ by attaching an external collar,
we can attach to $h'_t:H_0\to H_1$
a product map between the collars, getting a certain
$h''_t:H'_0\to H'_1$ whose trace on the boundary is $f_t$.

Since $f_t(\Gamma_0)=\Gamma_1$,
up to identifying $H'_i$ to $H_i$, the isotopy $\big(h''_t\big)_{t\in[0,1]}$
is now the desired one. To be completely formal and make sure that
$h''_j=h_j$ for $j=0,1$ one should at the beginning
slightly perturb $h_0$ and
choose internal collars
of $\partial H_i$ on which $h_0$ and $h_1$ behave as products, and then
identify $H'_i$ to $H_i$ just by rescaling the two-sided collar of
$\partial H_i$ in $H'_i$ to its one-sided collar in $H_i$.
\end{proof}

\paragraph{Hyperbolic volume}
Of course each element of $\calH_n$ (a relative handlebody of
the form $N(T)$ with $T\in \calT_n$) has volume
$n\cdot v_O$, where $v_O$ is the volume
of the hyperbolic regular ideal octahedron. Moreover $v_O=8\Lambda(\pi/4)\approx 3.66386$,
as first
computed by Milnor in~\cite{thurston:notes}.

If $T\in\calT_n$ then the numbers of $0$- and $1$-handles used to
construct $N(T)$ are respectively $n$ and $2n$, so $N(T)$ has
genus $n+1$ and its boundary has Euler characteristic $-2n$.
Miyamoto~\cite{Miy} has shown that $\calH_n$ is precisely the set
of hyperbolic manifolds $Y$ having minimal volume among those with
$\chi(\partial Y)=-2n$.

\begin{cor}
$\calH_n$ consists of the hyperbolic relative handlebodies of minimal volume among
those based on a handlebody of genus $n+1$.
\end{cor}

\paragraph{Complexity}
We now turn to the pant-meridinal complexity defined in the Introduction, establishing
all the remaining assertions of Theorem~\ref{T:to:H:intro:teo}.

\begin{wrapfigure}[13]{R}{4 cm}
    \begin{center}
    \mettifig{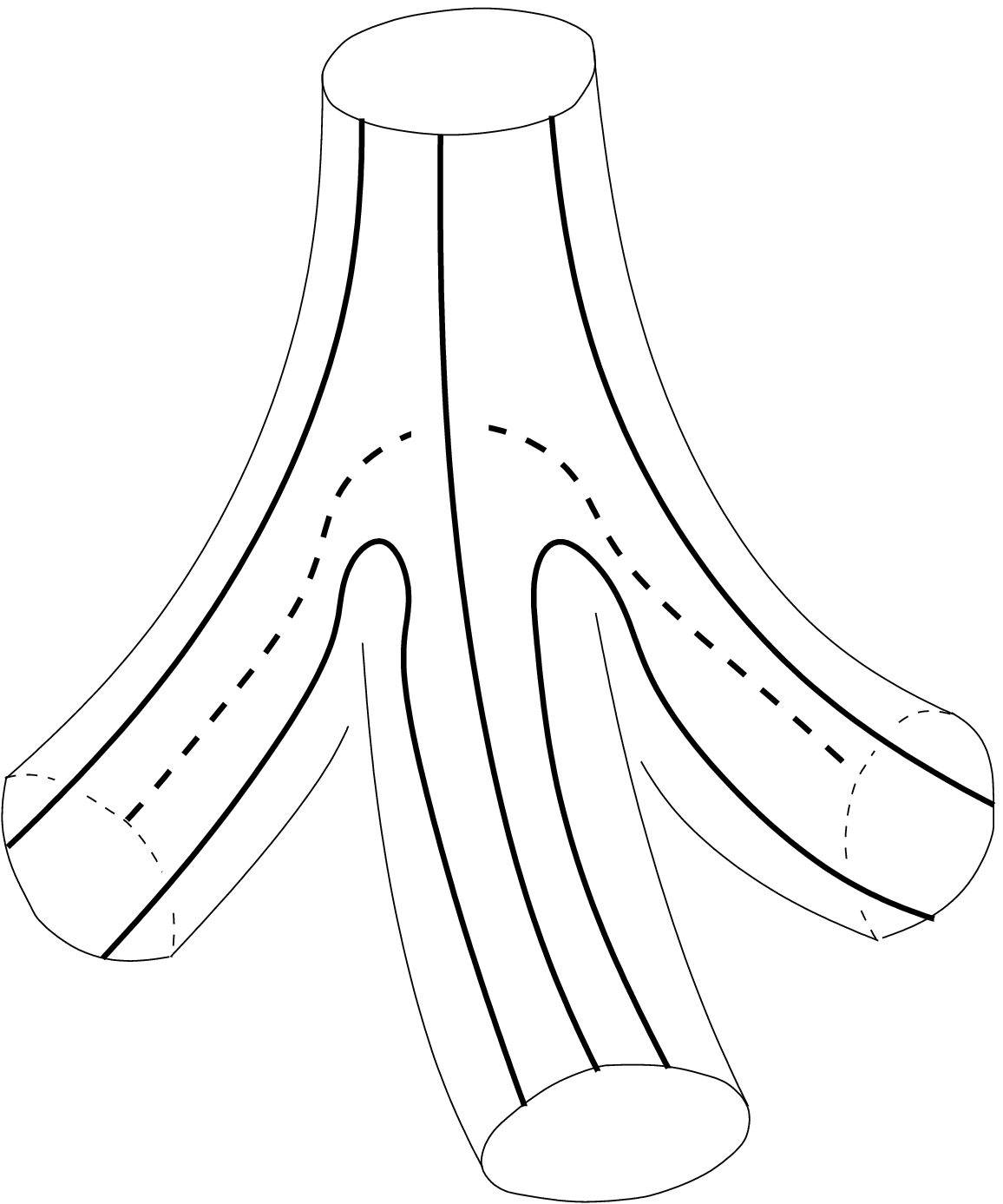,width=3cm}
    \nota{A tetrapod.} \label{frisia:fig}
    \end{center}
    \end{wrapfigure}

If $T\in\calT_n$ and $N(T)=(H,\Gamma)$, consider first the $2n$
meridinal discs of the $1$-handles used to construct $H$. These
discs cut $H$ into a union of ``tetrapods'' as shown in
Fig.~\ref{frisia:fig}.
Note that $\Gamma$ meets the boundary of each disc in three points.
To get a decomposition as in the definition of $c(\Gamma)$ we need to cut
each tetrapod along another disc. There are three non-isotopic but equivalent discs
one could choose, and each meets $\Gamma$ in $4$ points. This decomposition gives
the estimate $c(\Gamma)\leqslant 10n$.

We now claim that if $(H,\Gamma)$ is hyperbolic and $H$ has genus $n+1$
then $c(\Gamma)\geqslant 10n$.
Assume the contrary and pick a family of
discs which decomposes $H$ into
solid pairs of pants and meets $\Gamma$ in less than $10n$ points. Then there must exist one
solid pair of pants whose three meridinal
discs meet $\Gamma$ in $k<10$ points. But $k$ is even, so $k\leqslant 8<9=3\times 3$, whence
there is one disc meeting
$\Gamma$ in $h\leqslant 2$ points. Depending on whether $h$ is $0$, $1$, or $2$, this
disc gives either a boundary-compressing disc in $(H,\Gamma)$, or a disc compressing a toric (or Klein bottle)
cusp in
the genuine double of $(H,\Gamma)$, or an essential annulus in the double of $(H,\Gamma)$.
Since $(H,\Gamma)$ and its double are
hyperbolic, we get a contradiction.

We are left to show that a hyperbolic relative handlebody $(H,\Gamma)$ with $H$
of genus $n+1$ and $c(\Gamma)=10n$ arises as $N(T)$ for some
triangulation $T$. Let us consider the decomposition of $H$ into
solid pairs of pants which realizes the minimum $c(\Gamma)$. By
what already proved, for every pair of pants in the decomposition
there are two discs meeting $\Gamma$ in three points, and one disc
meeting $\Gamma$ in four points. It is now easy to show that there
is only one such configuration such that the number of
intersections between $\Gamma$ and the discs cannot be reduced,
whence the conclusion at once.

\section{Mirrored relative handlebodies}\label{link:section}
In this section we study the class $\calD$ of hyperbolic 3-manifolds
obtained by the mirroring construction applied to the elements of $\calH=N(\calT)$,
as described in the Introduction. The class $\calD$ was first defined in~\cite{Co:Th}.

\paragraph{Spaces associated to triangulations}
There is a common idea underlying the construction of the sets $\calH$ and $\calD$ of manifolds
in one-to-one correspondence with the set $\calT$ of triangulations.
It consists in choosing a topological or geometric block $B$ with four distinguished subsets
called \emph{faces}, so that $B$ has the same symmetries as
the tetrahedron, and in a procedure that
translates a combinatorial pairing between faces of tetrahedra 
into a homeomorphism between ``faces'' of the corresponding blocks $B$.
Using this procedure we associate to a triangulation $T$ a space $\Theta(T,B)$
obtained by gluing along faces some copies of $B$ according to the combinatorics of $T$.
For several choices of $B$ we can then prove that the topology or geometry of $\Theta(T,B)$
determines $T$ uniquely.

We describe now four examples of blocks $B$. 
\begin{itemize}
\item Take $\Delta$ itself as $B$, with tautological faces and gluing 
procedure. Then $\Theta(T,\Delta)$ is
the support $|T|$ of $T$;
    \begin{figure}
    \begin{center}
    \mettifig{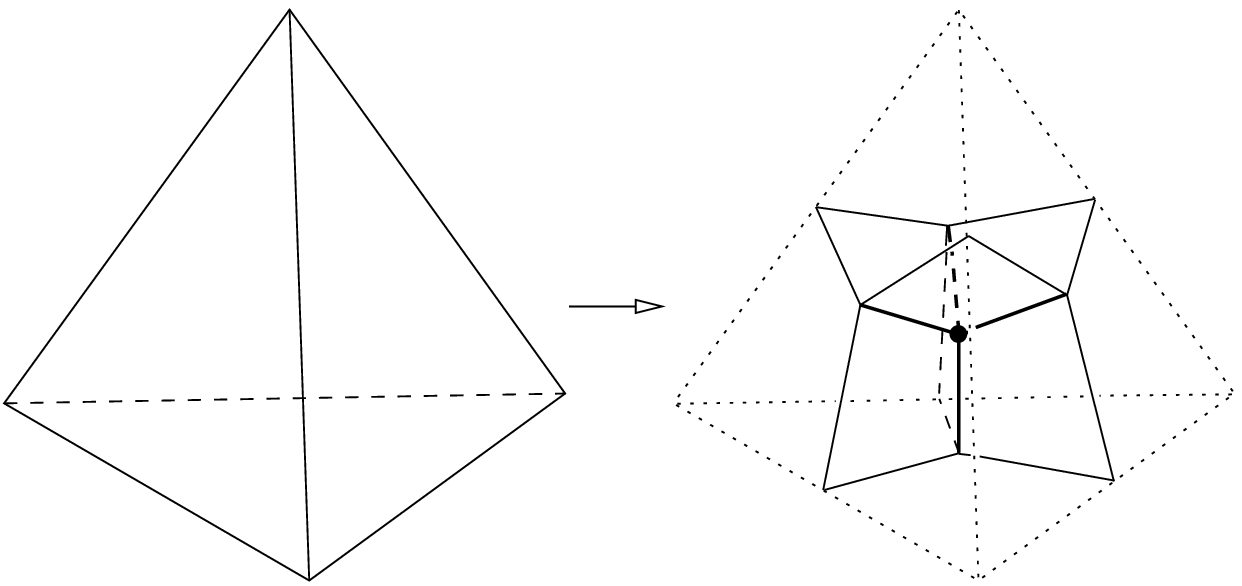,width=5 cm} \nota{From a triangulation to the dual special polyhedron.}
    \label{dualspine:fig}
    \end{center}
    \end{figure}
\item Within each of the four triangular faces of the tetrahedron $\Delta$
consider the \textsf{Y}-shaped
dual graph, and take $B$ as the cone from the centre of $\Delta$ over the union of these 
four \textsf{Y}-shaped graphs, which are defined to be the faces of $B$.
Face-pairings induce homeomorphisms between the corresponding \textsf{Y}-shaped graphs, and
the resulting space
$\Theta(T,B)$ is a \emph{special polyhedron} $P(T)$ called the \emph{dual} of $T$;
\item Take as $B$ the regular ideal octahedron $O$ as in Fig.~\ref{tetra:octa:fig}-right,
with white triangles as faces. Also here, face-pairings induce isometries between white triangles. 
The resulting metric space $\Theta(T,O)$ is then isometric to
$N(T)$, as shown in Section~\ref{tria:section};
\item Take as $B$ the block $K$ defined by mirroring $O$ in the
black triangles. This block was already considered in~\cite{Agol:K-blocks}, and it was described above
in Fig.~\ref{example:fig}.
The white triangles of $O$ glue up to $4$ geodesic thrice-punctured
spheres, which we define to be the faces of $K$. 
By Proposition~\ref{T:is:candeco:prop},
the isometry group of $K$ is
$$\{{\rm isometries\ of}\ O\ {\rm preserving\ colourings}\}\times\matZ/_{2\matZ}
\cong\permu_4\times\matZ/_{2\matZ}$$
where the non-trivial element of $\matZ/_{2\matZ}$
exchanges the two mirror copies of $O$.

We fix an arbitrary orientation on $B$ and note
that a face-pairing induces a correspondence between two thrice-punctured spheres, together with 
a bijection between their triples of punctures. We then glue the thrice-punctured spheres via the unique 
orientation-reversing isometry matching this bijection.
The resulting space
$\Theta(T,K)$ is oriented, and it is not hard to check that it is actually
isometric to the manifold $D(T)$, the orientable ``double'' of $N(T)$
defined in the Introduction.

\end{itemize}

It is a striking fact that, although $|T|$ does not determine $T$,
the topology of the other spaces $P(T)$, $N(T)$, and $D(T)$ does indeed
determine $T$. This is easy for $P(T)$, it was shown in Section~\ref{tria:section} for $N(T)$, and
it is proved below for $D(T)$:

\begin{teo}\label{D:bijective:teo}
The map $D:\calT\to\calD$ is bijective.
\end{teo}

\begin{figure}
\begin{minipage}{.03\textwidth}\hfil
\end{minipage}
\begin{minipage}{.3\textwidth}
\centering
\mettifig{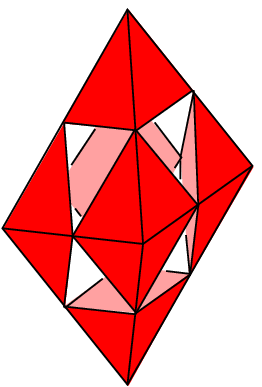, width = .9\textwidth}
\vglue -4mm
\end{minipage}
\begin{minipage}{.03\textwidth}\hfil
\end{minipage}
\begin{minipage}{.57\textwidth}
{The decomposition of $D(T)$ into blocks $K$, and hence into octahedra,
gives a horospherical cusp neighbourhood of $D(T)$, constructed by taking in each octahedron
the horospherical neighbourhoods of the vertices shown here.
Each such neighbourhood has volume equal to $1/2$ and is based on a unitary Euclidean square.}
\end{minipage}
\nota{The horospherical cusp neighbourhood of $D(T)$ given by $T$.}
\label{cusp_section2:fig}
\end{figure}

As a preliminary to the proof of this result,
note that $D(T)$ is a cusped orientable finite-volume hyperbolic $3$-manifold without boundary.
Its volume is $2v_O$ times the number of
tetrahedra of $T$, \emph{i.e.}~twice the volume of $N(T)$.
The cusps of $D(T)$ and $N(T)$ correspond to the edges of $T$, and
the geometry of a cusp of $D(T)$ is determined
by the geometry of the corresponding cusp in $N(T)$, which is
a Euclidean annulus or M\"obius band.
Now recall that in general,
if $V$ is an orientable hyperbolic $3$-manifold with cusps,
a \emph{horospherical cusp neighbourhood} $\calO$ of $V$ is a union of
disjoint open subsets of $V$, one in each cusp, which lift to open horoballs in $\matH^3$.
Note that $\partial\calO$ is a family of
immersed Euclidean tori which intersect each other and self-intersect at most tangentially.
The following result is an easy consequence of Lemma~\ref{strip:shape:lem}.

\begin{figure}
\begin{center}
\mettifig{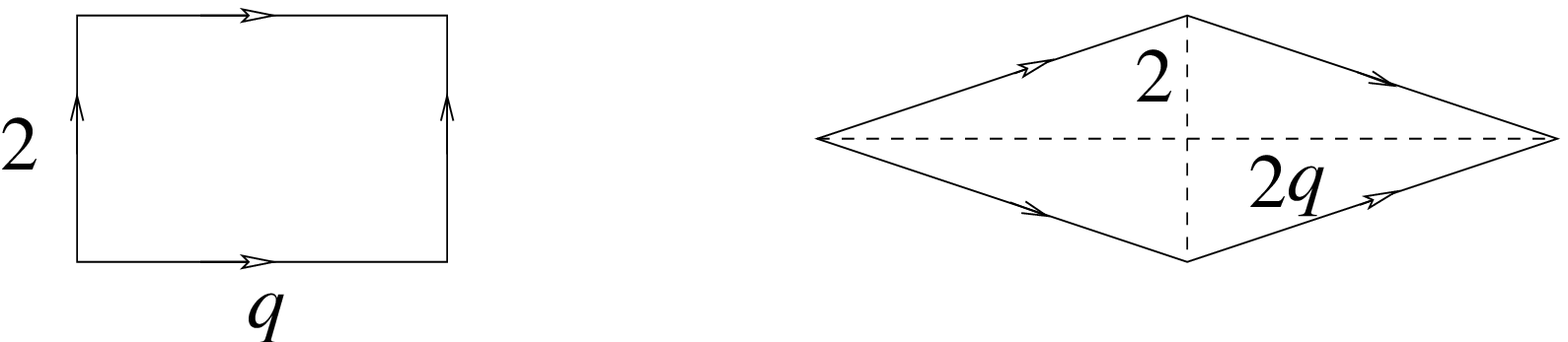,width=8 cm}
\nota{Shape of a boundary component of the horospherical cusp neighbourhood of $D(T)$. The
rectangle occurs if the corresponding cusp of $N(T)$ is orientable, the rhombus if it is
non-orientable.} \label{cusps:fig}
\end{center}
\end{figure}

\begin{cor}\label{geometry:of:cusps:cor}
Let $T\in\calT$ and let $\calO(T)$ be the horospherical cusp neighbourhood defined in
Fig.~\ref{cusp_section2:fig}. Then $\calO(T)$ is maximal. Moreover, if a component $E$ of
$\partial\calO(T)$ bounds a cusp which corresponds to an edge of $T$ with valence $q$, then $E$
is a Euclidean torus of area $2q$ and geometric
shape as described in Fig.~\ref{cusps:fig}.
\end{cor}

\begin{figure}
\begin{center}
\mettifig{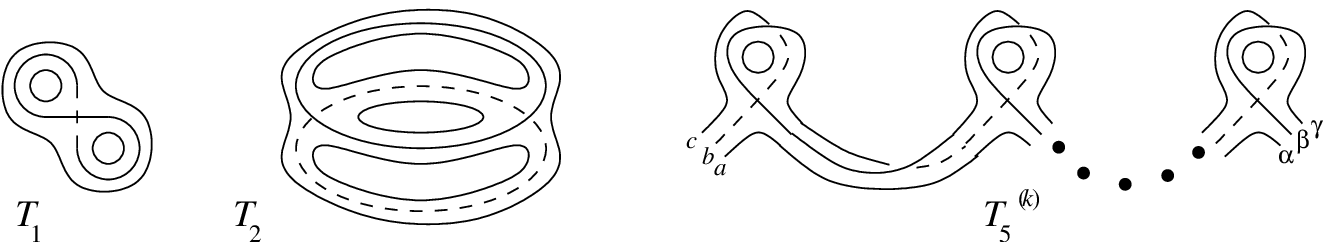,width=12cm}
\nota{Particular triangulations. The picture
describes the dual special polyhedron by means of the attaching circles of the
regions to the singular set.} \label{combinations:fig}
\end{center}
\end{figure}

For the proof of Theorem~\ref{D:bijective:teo} and later in Section~\ref{isom:section}
we will need to refer to certain specific triangulations $T_1$,
$T_2$, $T^{(k)}_3$, and $T^{(k)}_4$ for $k\geqslant 1$. The first two are defined in
Fig.~\ref{combinations:fig}-left and centre, while $T^{(k)}_3$ (respectively, $T^{(k)}_4$)
is obtained from the triangulation $T^{(k)}_5$ of
Fig.~\ref{combinations:fig}-right by matching $(a,b,c)$ to $(\alpha,\beta,\gamma)$ (respectively,
$(\alpha,\gamma,\beta)$). Both $T^{(k)}_3$ and $T^{(k)}_4$ contain $k$ tetrahedra. As an easy
application of Corollary~\ref{geometry:of:cusps:cor} one can now show the following:

\begin{cor}\label{T8T9:cor}
For $j=3,4$ and $k\geqslant 1$ the manifolds $D(T^{(k)}_j)$ are pairwise non-homeomorphic.
\end{cor}
\begin{wrapfigure}[9]{r}{4cm}
\vspace{-.5 cm}
\begin{center}
\mettifig{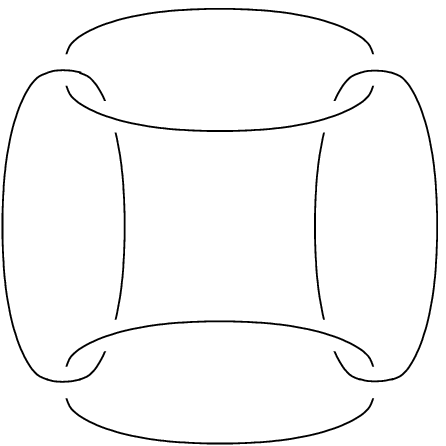,width=3cm}
\nota{A link in $S^3$} 
\label{chainfour:fig}
\end{center}
\end{wrapfigure}

We also note that $D(T_1)$ is a rather remarkable
manifold. For instance, one can see that it is the complement
in $S^3$ of the link shown in Fig.~\ref{chainfour:fig}, 
and that its isometry group, which has 64 elements, 
acts as the dihedral group $D_4$ on the cusps.

The proof of Theorem~\ref{D:bijective:teo} is an immediate consequence of the following:

\begin{prop}\label{D:bijective:core:prop}
Let $V$ be an orientable cusped hyperbolic $3$-manifold.
\begin{enumerate}
\item Every decomposition of $V$ into blocks $K$ arises from a unique triangulation $T$ such that $V=D(T)$;
\item If $V\not\in\{D(T_1),D(T_2),D(T^{(k)}_3),D(T^{(k)}_4):\ k\geqslant 1\}$ then the decomposition of $V$ into blocks $K$ is unique;
\item If $V\in\{D(T_1),D(T_2),D(T^{(k)}_3),D(T^{(k)}_4):\ k\geqslant 1\}$ then $V$ has
distinct decompositions into blocks $K$, but they arise from combinatorially equivalent triangulations.
\end{enumerate}
\end{prop}

\dimo{D:bijective:core:prop}
The first assertion is easy: in a decomposition of $V$ into
blocks $K$ all the gluings are orientation-reversing isometries, and
the set of these isometries corresponds bijectively to $\permu_4\cong\Aut(\Delta)$.
Therefore the procedure described above to pass from $T$ to $\Theta(T,K)=D(T)$ can be reversed,
and we are done.

We establish assertions 2 and 3 at the same time, the proof of
the latter being an easy by-product of the
proof of the former. We give two different arguments, one analyzing the positions
of geodesic thrice-punctured spheres in $V$, and one using the Epstein-Penner decomposition.

For the first proof, assume that $V$ has two distinct
decompositions $\rho$ and $\rho'$ into blocks $K$, arising from triangulations
$T$ and $T'$. We must show that $T$ and $T'$ are isomorphic to
each other and to one of $T_1,T_2,T^{(k)}_3,T^{(k)}_4$, $k\geqslant 1$.
Note first that $T$ and $T'$ must have the same number $n$ of tetrahedra.
For the sake of simplicity we first
assume $n>1$. Later we will sketch how to treat the easy case of triangulations with
a single tetrahedron.

Remark now that the decompositions $\rho$ and $\rho'$ are made along $2n$ geodesic
thrice-punctured spheres, and let $S_1,\ldots,S_{2n}$ be those giving $\rho$.
Since $\rho'\neq\rho$ there is another such sphere $S'$ distinct from the $S_i$'s.
Our first aim is to analyze how $S'$ can intersect the $S_i$'s.
Before proceeding, note that the thrice-punctured sphere has a unique hyperbolic structure,
and it contains precisely $6$ complete simple
geodesics.

\begin{figure}
\begin{center}
\mettifig{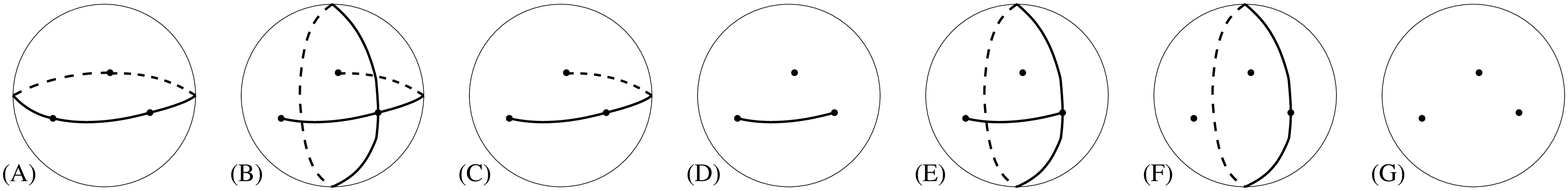,width=12 cm}
\nota{Disjoint systems of geodesics in the thrice-punctured sphere.} \label{spheres:fig}
\end{center}
\end{figure}

\begin{figure}
\begin{center}
\mettifig{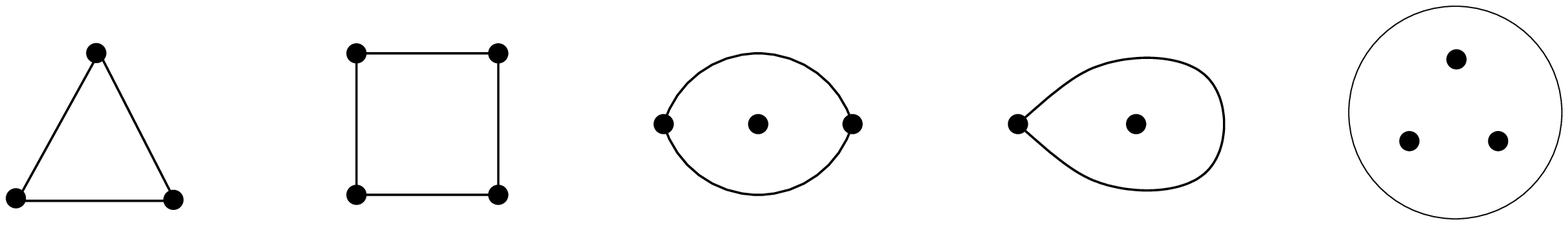,width=8 cm}
\nota{Surfaces obtained cutting $S'$ along $\lambda$.} \label{spheres2:fig}
\end{center}
\end{figure}

\begin{figure}
\begin{center}
\mettifig{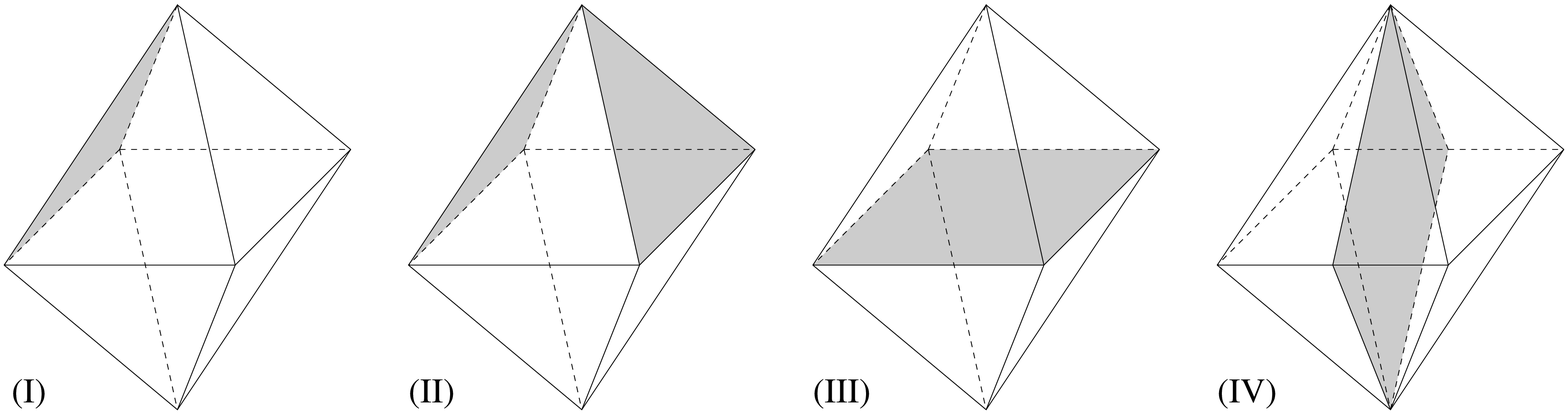,width=10 cm}
\nota{The intersection between $S'$ and one of the blocks $K$. In all four cases we are only showing the
octahedron $O$ of which $K$ is the double. In cases (I)-(III) the intersection
$S'\cap K$ is contained in $O$, while in case (IV) it contains another
similar quadrilateral in the mirror copy of $O$.} \label{poligonigeod:fig}
\end{center}
\end{figure}

Since $S'$ and $\bigcup S_i$ are geodesic surfaces, they intersect transversely in a
disjoint union $\lambda$ of geodesic lines.
Therefore $\lambda$ as a subset of $S'$
gives one of the configurations (A)-(G) shown in Fig.~\ref{spheres:fig}.
Cutting $S'$ along $\lambda$ we then get one or two of the surfaces shown in Fig.~\ref{spheres2:fig},
and each such surface must be contained in one of the blocks $K$
arising from $T$. We now note that $\lambda$ is a union of
geodesics also in $\bigcup S_i$, and $S_i$ is the double of one of the
white faces of the octahedron $O$. It readily follows that the
intersection between $S'$ and each white triangle is either an
edge or a height (the bisecant of a corner) of the triangle.
Combining this fact with the remark that $S'\setminus\lambda$
contains at most two components as in Fig.~\ref{spheres2:fig}, we
see that the intersection between $S'$ and each of
the blocks $K$ is empty or as in Fig.~\ref{poligonigeod:fig}.

The components of $S'\cap K$ which appear in Fig.~\ref{poligonigeod:fig}
are either ideal squares, or ideal triangles, or punctured ideal bigons.
It follows that, among the configurations of $\lambda$ on $S'$ shown in Fig.~\ref{spheres:fig},
only cases (A)-(D) are possible. Moreover case (A) can be realized either by two triangles
of type (I) in different blocks $K$, or by a pair of triangles of type (II) in one block.
We call these configurations (A$'$) and (A$''$) respectively. In a similar way we have
cases (B$'$) and (B$''$), while case (C) can only arise from one square of type (III),
and case (D) only from a pair of quadrilaterals of type (IV).

\begin{figure}
\begin{center}
\mettifig{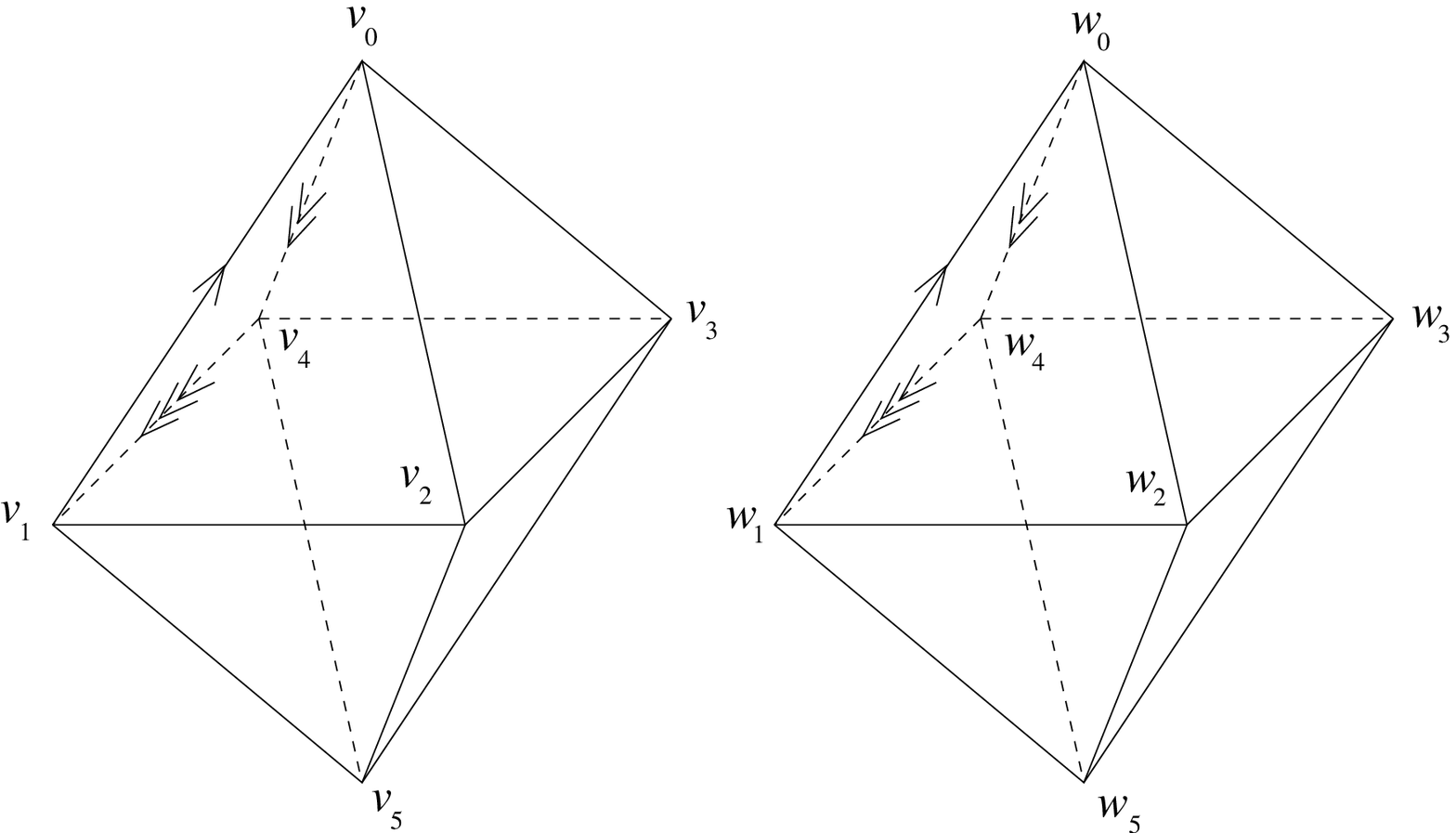,width=8 cm}
\nota{Configuration (A$'$).} \label{caseA:fig}
\end{center}
\end{figure}

\begin{figure}
\begin{center}
\mettifig{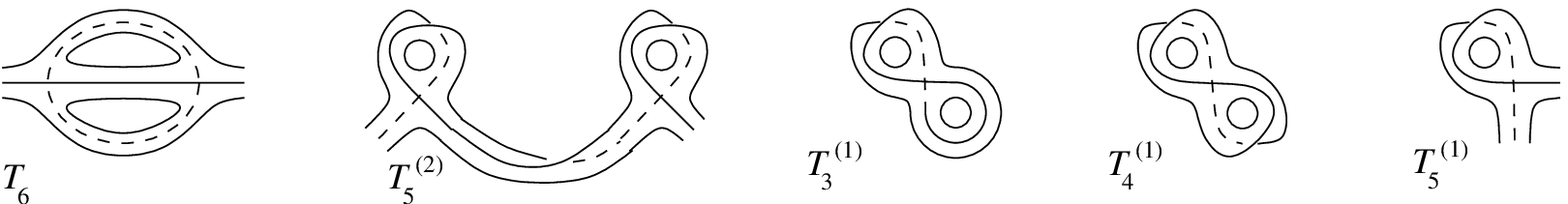,width=12.5cm}
\nota{Particular triangulations or fragments of.} \label{fragments:fig}
\end{center}
\end{figure}

We now prove that the existence
of a thrice-punctured sphere $S'$ of one of the types (A$'$), (A$''$),
(B$'$), (B$''$), (C), or (D) forces a portion of $T$ (or the whole of it) to have
some definite shape. For instance, in case (A$'$) there must be two blocks $K$
with gluings between their boundary thrice-punctured spheres that induce the
edge-identifications of Fig.~\ref{caseA:fig}.
Since gluings are orientation-reversing,
the edge marked with a single arrow forces the
gluing $(v_0,v_1,v_2)\to(w'_0,w'_1,w'_2)$, where
$(v_0,v_1,v_2)$ is a white triangle of the octahedron on the left, and
$(w'_0,w'_1,w'_2)$ is the mirror image of $(w_0,w_1,w_2)$ in the companion of the octahedron
on the right. Similarly we must have gluings
$(v_0,v_3,v_4)\to(w'_0,w'_3,w'_4)$ and $(v_5,v_1,v_4)\to(w'_5,w'_1,w'_4)$.
It easily follows that $T$ contains a subtriangulation $T_6$
as in Fig.~\ref{fragments:fig}.
A long but straight-forward analysis shows that case (A$''$) leads
to $T_1$ (defined above in Fig.~\ref{combinations:fig}), case
(B$'$) to $T_5^{(2)}$, case (B$''$) to $T_3^{(1)}$ or $T_4^{(1)}$, case (C) to $T_1$
again, and case (D) to $T_5^{(1)}$. Cases $T_1$, $T_3^{(1)}$, and $T_4^{(1)}$ are
actually forbidden by the assumption $n>1$, so we discard them for
the time being.

\begin{figure}
\begin{center}
\mettifig{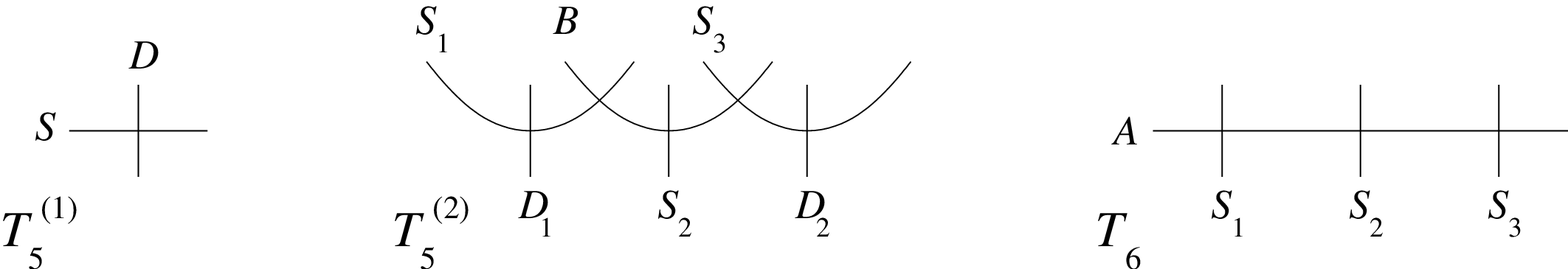,width=11cm}
\nota{Patterns of intersections of surfaces in $\calS$.} \label{intersections1:fig}
\end{center}
\end{figure}

Having analyzed completely how a single geodesic thrice-punctured sphere $S'$ can appear
with respect to the $S_i$'s, we now consider the collection $\calS$ of all such spheres (including the $S_i$'s).
We begin by noting that the portions of triangulation found so far
contain several elements of $\calS$, namely
$T_1$ contains two surfaces of type (A$''$)'s, four of type (C), and two $S_i$'s,
$T_3^{(1)}$ contains one  (B$''$), one (D) and two $S_i$'s,
$T_4^{(1)}$ contains two  (B$''$)'s, two (D)'s and two $S_i$'s,
$T_5^{(1)}$ contains one (D) and one $S_i$,
$T_5^{(2)}$ contains one (B$'$), two (D)'s and three $S_i$'s,
and $T_6$ contains one (A$'$) and three $S_i$'s.
This implies first of all that $\calS$ is finite.
Moreover we can study how the various elements of $\calS$ contained in the fragment
intersect each other, which always occurs along one or more geodesics. The resulting
configurations are shown by the `incidence graphs' of Fig.~\ref{intersections1:fig} (where we omit
$T_1$, $T_3^{(1)}$, and $T_4^{(1)}$, which are forbidden at this stage).

\begin{figure}
\begin{center}
\mettifig{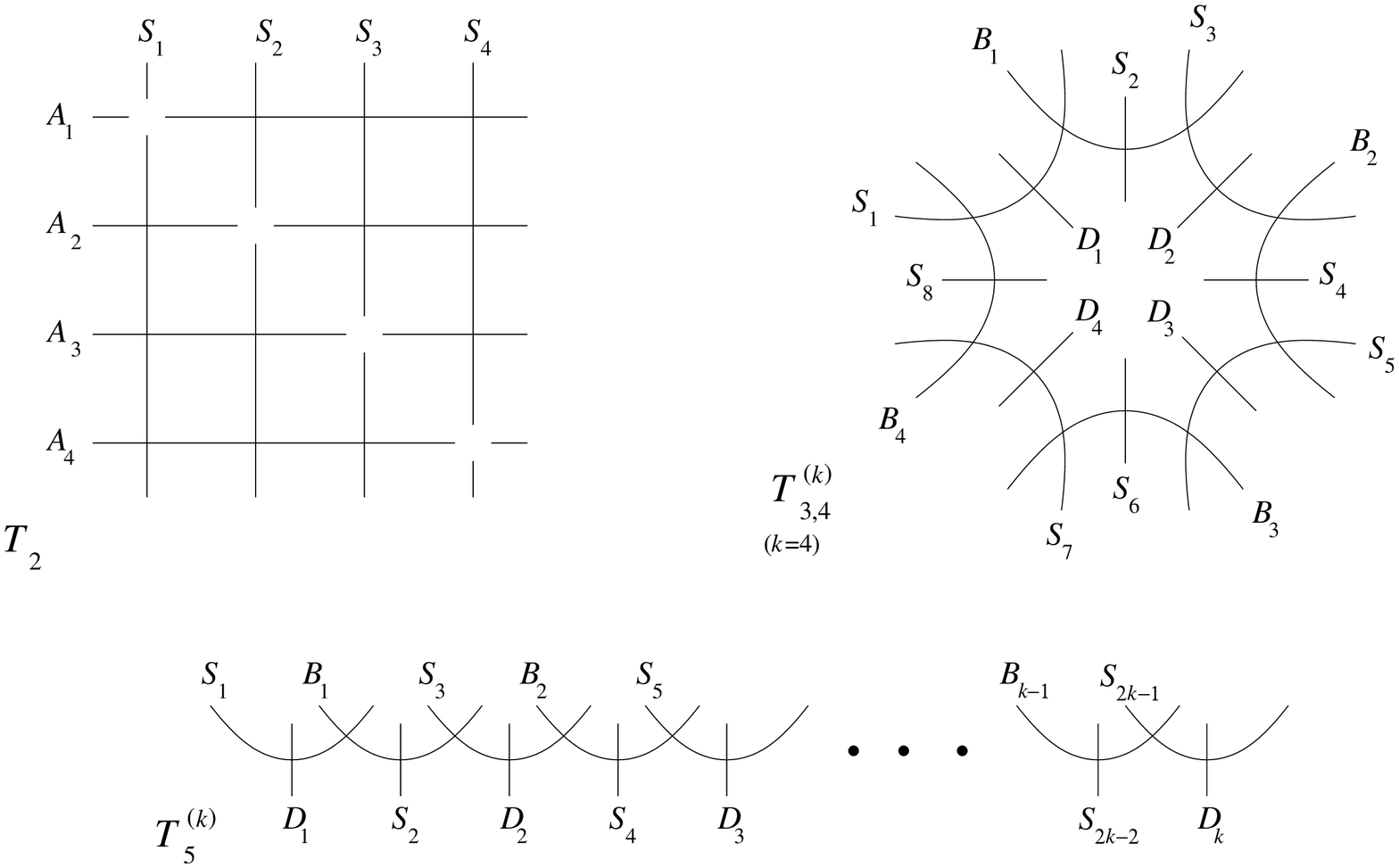,width=12.5cm}
\nota{Patterns of intersections of elements of $\calS$ in $T_2,T^{(k)}_3,T^{(k)}_4,T^{(k)}_5$.}\label{intersections2:fig}
\end{center}
\end{figure}

Let us consider now a connected component of $\bigcup\calS$. By what shown already this component
must be contained in a portion of $T$ obtained by
assembling fragments as in Fig.~\ref{fragments:fig}. This portion of $T$ is then one of
$T_2,T^{(k)}_3,T^{(k)}_4$ with $k\geqslant 2$, or $T^{(k)}_5$ with $k\geqslant 1$, or $T_6$.
The incidence graphs of the elements of $\calS$ contained in these triangulations are shown in
Fig.~\ref{intersections2:fig} (except $T_6$, already shown above).

Let us note now that the incidence graphs which appear in Figg.~\ref{intersections1:fig}
and~\ref{intersections2:fig} and can occur in $V$ are pairwise distinct
(remark that $T^{(k)}_3$ and $T^{(k)}_4$
have the same incidence graph, but $V$ cannot be $D(T^{(k)}_3)$ and $D(T^{(k)}_4)$
at the same time, by Corollary~\ref{T8T9:cor}).
Recall that we are assuming that $V$ has another decomposition $\rho'$ into blocks $K$.
The above discussion on how the set of all thrice-punctured spheres intersects the spheres
giving $\rho$ applies also to the spheres giving $\rho'$. We deduce that in the intersection
graph of some component of $\calS$ there should exist an automorphism not leaving the
$S_i$'s invariant. This can only happen in cases $T_2$, $T^{(k)}_3$, $T^{(k)}_4$, and $T_5^{(1)}$.
However in $T_5^{(1)}$ the r\^ole of $S_i$ can be identified intrinsically,
because within it $S_i\cap D$ joins a
puncture of $S_i$ to itself, while it joins two distinct punctures of $D$.
This eventually establishes assertion 2 for manifolds not decomposing in a single block $K$.
To prove assertion 3 we only need to note that, as one sees from the incidence graphs,
$T_2$, $T^{(k)}_3$, and $T^{(k)}_4$ contain precisely
two families of punctured spheres giving a decomposition into blocks $K$, and these
families arise from combinatorially equivalent triangulations.

We are left to deal with the case where $V=D(T)$ and $T$ has one tetrahedron.
As a by-product of the above argument we see that the decomposition into blocks $K$
could be non-unique at most if $T$ is $T_1$, $T_3^{(1)}$, or $T_4^{(1)}$. Now it is clear
that $D(T_1)$, $D(T_3^{(1)})$, and $D(T_4^{(1)})$ are pairwise distinct, because they contain
respectively 8, 6, and 4 thrice-punctured spheres. So at least the conclusion
that $D(T)$ determines $T$ is clear. A direct argument, again based on the analysis of the
intersection graphs of the family $\calS$, proves that in all these cases multiple
decompositions indeed exist.

As mentioned in the Introduction, we have an alternative and more geometric proof
of the proposition, which we believe is worth explaining. This proof is based on the polyhedral
decomposition of $V=D(T)$ due to Epstein
and Penner~\cite{EpPe}, which we now briefly recall.
Let $\calO$ be a horospherical cusp neighbourhood for $V$.
If $\mathrm{Cut}(V,\calO)$ is the cut-locus of $V\setminus\calO$
relative to its boundary, the Epstein-Penner decomposition of $V$ relative to $\calO$
is given by the ideal polyhedral decomposition of $V$
dual to $\mathrm{Cut}(V,\calO)$.
Let $\calO(T)$ be as defined in Fig.~\ref{cusp_section2:fig}.
Let $\delta(T)$ denote the geometric decomposition of $V$
into ideal octahedra induced by $T$.
It is easy to see that $\delta(T)$ is the Epstein-Penner decomposition
of $V$ relative to $\calO(T)$.

Let $\{e_i\}$ be
the edges of $T$, and $\{q_i\}$ be the corresponding valences.
Let $C_i$ be the toric cusp of $V$ corresponding
to $e_i$. The component of $\calO(T)$ corresponding to
$C_i$ has volume $q_i$. Observe that
the total volume of
$\calO(T)$ is equal to $3\cdot\#\delta(T)=3\cdot\mathrm{Vol}(V)/v_O$,
so it depends solely on $V$.

Now let us suppose $V=D(T)=D(T')$ with $T\neq T'$.
Let $e'_i$ be the edge of $T'$ corresponding
to $C_i$ and let $q'_i$ be the multiplicity of $e'_i$ in $T'$.
If $q_i=q'_i$ for all $i$ then
$\calO(T)=\calO(T')$, so $\delta(T)=\delta(T')$ since
$\delta(T)$ and $\delta(T')$ both give the Epstein-Penner decomposition
of $V$ relative to $\calO(T)=\calO(T')$.
Moreover, the octahedra of $\delta(T)$ glue up in pairs into blocks $K$
in a unique way, except when $T$ is the triangulation $T_2$:
two decompositions arise in this case, both induced by the same $T$ anyway.

\begin{figure}
\begin{center}
\mettifig{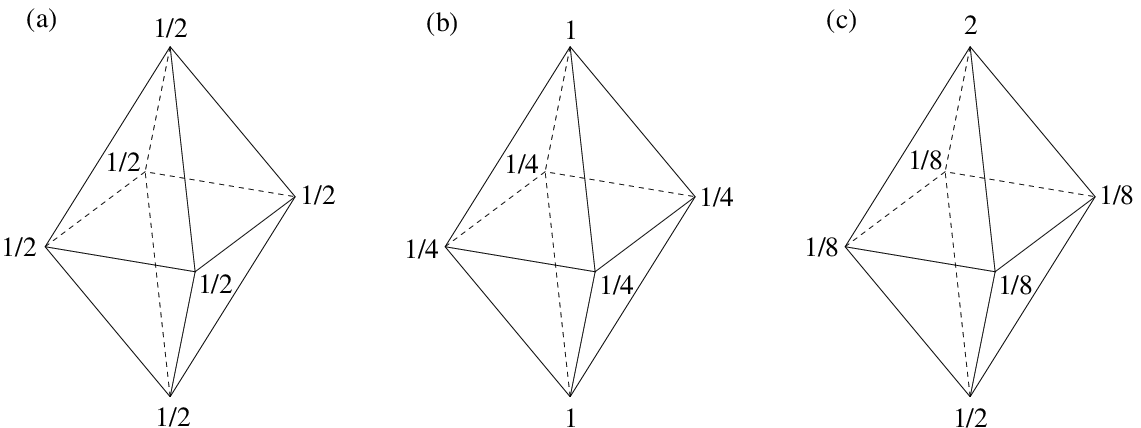}
\nota{The allowed values for $r$ on the vertices of an octahedron
in $T$.} \label{octaallow:fig}
\end{center}
\end{figure}

So let us suppose $q_i\neq q'_i$ for some $i$.
Corollary~\ref{geometry:of:cusps:cor}
implies that one of the following conditions holds:
\begin{itemize}
\item $\{q_i, q'_i\}=\{1,4\}$ and the
cusps of $N(T)$ and $N(T')$ along $e_i$ and $e'_i$ are both orientable;
\item $q_i=2$, the cusp of
$N(T)$ along $e_i$ is orientable, $q'_i=1$, and the cusp of
$N(T')$ along $e'_i$ is non-orientable; or viceversa.
\end{itemize}
Let $v$ be a vertex of an octahedron $O$ in $\delta(T)$.
We define $r(v)$ to be the
volume of the connected component
of $\calO(T')\cap O$ which is asymptotic to $v$.
The previous observation implies that
$r(v)\in\{1/8,1/4,1/2,1,2\}$.
Note now that the components of $\calO(T')$ lying within a certain
octahedron $O\in \delta(T)$ do not overlap, and a direct computation shows
that the volume of $\calO(T')\cap O$ is then at most $3$. Moreover the total volume of
$\calO(T')$ is thrice the number of octahedra of $\delta(T)$, so the volume of
$\calO(T')\cap O$ is precisely $3$.
This easily implies
that every octahedron in $\delta(T)$ belongs to one of
the three classes described in Fig.~\ref{octaallow:fig},
according to the values that $r$ takes on its vertices.
If two vertices of octahedra of $T$
are glued to each other then $r$ takes on them the same value.
Since $V$ is connected we deduce that all the octahedra
of $T$ are of the same type.

If all the tetrahedra are of type (a)
then $q'_i=q_i$ for all $i$, whence $T=T'$, as seen above.
To analyze cases (b) and (c), note first
that a vertex $v$ with $r(v)=2$ gives an edge of valence $1$ in $T$,
\emph{i.e.}~only one octahedron of $\delta(T)$ is incident to $v$, and
the corresponding cusp of $N(T)$ is orientable.
Similarly if $r(v)=1/8$ then $v$ gives an ``orientable" edge of valence $4$, if $r(v)=1$ then
$v$ gives a ``non-orientable" edge of valence $1$,
if $r(v)=1/4$ then
$v$ gives an ``orientable" edge of valence $2$.
Using these facts it is easy to see that if there is in $\delta(T)$ an octahedron of type (b) then there is only one,
$T$ is $T_1$, and $T'$ is also $T_1$. Suppose then that all octahedra are of type (c).
If there are $k$ of them then $T$ and $T'$ are either
$T^{(k)}_3$ or $T^{(k)}_4$, and Corollary~\ref{T8T9:cor} implies that $T=T'$ also in this case.
\finedimo

Proposition~\ref{D:bijective:core:prop} enables
us to compute the isometry group of $D(T)$ in terms of the
group Aut$(T)$ of combinatorial automorphisms of $T$.

\begin{cor} \label{isometrie:cor}
Let $T$ be a triangulation distinct from $T_1$, $T_2$, $T^{(k)}_3$,and $T^{(k)}_4$. Then
$${\rm Isom}(D(T))\cong {\rm Isom}^+(D(T))\times \matZ/_{2\matZ},$$
$${\rm Isom}^+(D(T))\cong {\rm Isom}(N(T)) \cong {\rm Aut}(T).$$
\end{cor}
\begin{proof}
The isomorphism
${\rm Isom}(N(T)) \cong {\rm Aut}(T)$
is an immediate consequence of
Proposition~\ref{T:is:candeco:prop}.
By construction
every automorphism $\varphi\in\mathrm{Aut}(T)$ induces an
orientation-preserving isometry $\overline{\varphi}\in\mathrm{Isom}(D(T))$.
Moreover we can define an orientation-reversing
isometry $\tau$ of $D(T)$ by reflecting any block $K$ along
the black faces of the octahedra the block is made of,
and of course $\tau\circ\overline{\varphi}=\overline{\varphi}\circ\tau$
for any $\varphi\in\mathrm{Aut}(T)$.
As shown in Proposition~\ref{D:bijective:core:prop},
if $T$ is not $T_1, T_2,T^{(k)}_3,T^{(k)}_4$,
then $D(T)$ admits a unique decomposition into blocks $K$, which is
therefore preserved by any element in $\mathrm{Isom}(D(T))$.
This implies that any isometry of $D(T)$ is induced, up to composition
with $\tau$, by an element
in $\mathrm{Aut}(T)$, and our assertions easily follow.
\end{proof}

\section{Dehn filling}\label{surg:section}
In this section we
continue our study of the class $\cal D$ of complete hyperbolic
$3$-manifolds with cusps, introduced in~\cite{Co:Th} and shown in
Section~\ref{link:section} to
be in one-to-one correspondence with the set $\calT$ of triangulations.

\paragraph{Meridians and longitudes}
Let $T$ be an arbitrary triangulation. To provide some information on the
Dehn fillings of $D(T)$,
we first recall from Corollary~\ref{geometry:of:cusps:cor}
that $D(T)$ has a maximal horospherical cusp
neighbourhood, bounded by a Euclidean torus of area $2q$ on each cusp.
Here, $q$ is the valence of the corresponding edge $e$ of $T$, and
the cusp has one of the two shapes shown in
Fig.~\ref{cusps:fig}. The
rectangle occurs if and only if the corresponding cusp in $N(T)$ is
orientable, or equivalently if the regular neighbourhood of $e$ in $M(T)$
is a manifold. Therefore, if $T$ is an ideal triangulation of a
manifold, only rectangles occur.

\begin{wrapfigure}[13]{r}{5 cm}
\vspace{-.4 cm}
\begin{center}
\mettifig{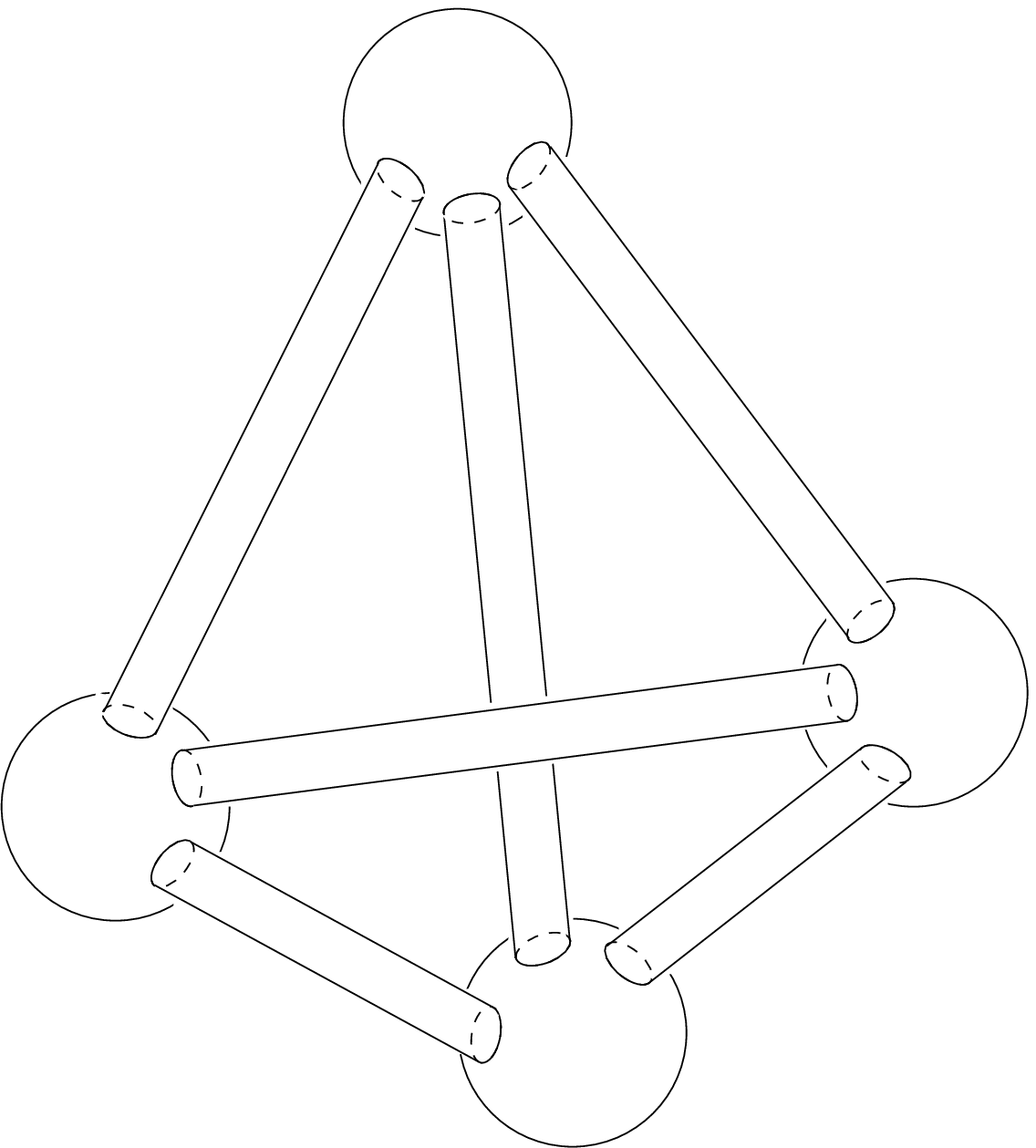,width=4cm} 
\nota{A handlebody.} 
\label{tubesandpants:fig}
\end{center}
\end{wrapfigure}

Recall that $T$ induces a decomposition of $D(T)$ into blocks $K$.
Each such block can be seen as the complement in $S^3$ of
the handlebody shown in Fig.~\ref{tubesandpants:fig} (the tubes give the cusps of $K$, and
the pairs of pants give the geodesic thrice-punctured spheres in the 
boundary of $K$) and it contributes
to the bases of the cusps with six Euclidean tubes of width 1 and boundary 
loops of length 2.

We define any boundary loop of a tube to be the \emph{meridian} of the
corresponding cusp.  In Fig.~\ref{cusps:fig} the meridians are given by the vertical
edges of the rectangle and by the vertical diagonal of the rhombus.
If $T$ is an ideal triangulation of a manifold,
a preferred \emph{longitude} on each cusp can also be defined as the unique
loop having length $q$, \emph{i.e.}~as the
loop given by the horizontal edges of the rectangle in Fig.~\ref{cusps:fig}-left.

Let $T$ consist of $n$ tetrahedra.
It is clear from Fig.~\ref{tubesandpants:fig} that by performing a
meridinal Dehn filling on each cusp we get a connected sum of
$n+1$ copies of $S^2\times S^1$, as also
stated in Proposition~\ref{X:prop}.
Therefore the first Betti number of $D(T)$ is at least $n+1$. More
precisely, we have the following:

\begin{prop} \label{homology:prop}
Let $T\in\calT$ consist of $n$ tetrahedra. Then
$$H_1(D(T);\matZ)=\matZ^{n+1}\oplus H^2(P(T);\matZ).$$
\end{prop}

\begin{proof}
Set $D=D(T)$ and $P=P(T)$. Denote by $S(P)$ the singular set of
$P$, \emph{i.e.}~the 4-valent graph representing the gluings in
$T$. Edges and vertices of $S(P)$ correspond to decomposing pairs
of pants and blocks $K$ of $D$. Let $p:D\to S(P)$ be the
composition of the obvious projections $D\to N(T)\to P_0(T)\to
S(P)$, and note that $p$ sends a neighborhood of each decomposing
pair of pants to the corresponding edge and each block $K$ to the
corresponding vertex. This gives a surjection
$p_*:H_1(D)\twoheadrightarrow H_1(S(P))=\matZ^{n+1}$. Therefore
$H_1(D) = \matZ^{n+1}\oplus\Ker( p_*)$.

Every loop in $\Ker(p_*)$ is homologous to a sum of loops each contained in
a block $K$. The cores of the six tubes of Fig.~\ref{tubesandpants:fig}
generate $H_1(K)$, so the meridians $\mu_i$, with arbitrary orientation, generate $\Ker( p_*)$.
Each pair of pants separating two blocks $K$
gives a relation of type $\pm\mu_{j_1}\pm\mu_{j_2}\pm\mu_{j_3}=0$, with appropriate signs.
We now claim that every other relation is superfluous.

A relation between the $\mu_i$'s comes from an embedded orientable surface bounded by
some copies of the $\mu_i$'s. This surface intersects each pair of pants
in a union of loops. The trivial loops can be dismissed by elementary cut-and-paste operations,
so we can assume each loop is parallel to some $\mu_i$. Therefore the relation can be expressed
as a sum of relations coming from a set of surfaces each contained in a block $K$. Now consider $K$
as the complement in $S^3$ of four balls $B_1,\ldots,B_4$
and six solid tubes, as in Fig.~\ref{tubesandpants:fig}. Adding discs to
a surface contained in $K$ we get a closed surface in $S^3\setminus(B_1\cup\ldots\cup B_4)$.
Any such surface compresses in $S^3\setminus(B_1\cup\ldots\cup B_4)$ to a union of spheres parallel to the $\partial B_i$'s,
so we can express the relation coming from a surface contained in $K$ as a sum
of relations induced by the pairs of pants in $\partial K$. This proves our claim.

We have shown that $\Ker(p_*)$ has one generator $\mu_i$ for each cusp and one relation
$\pm\mu_{j_1}\pm\mu_{j_2}\pm\mu_{j_3}=0$ for
each pair of pants. Cusps and pairs of pants correspond to faces and edges of $P$, and this gives
$\Ker(p_*)= H^2(P;\matZ)$.
\end{proof}

\paragraph{Dehn filling}
Let us fix a triangulation $T$ with $n$ tetrahedra.
By a Dehn filling of $D(T)$ we mean the result of some Dehn filling
on some (possibly all) cusps of $D(T)$. We begin with the following observation, needed below:

\begin{rem}\label{DM:rem}
{\em
If $T$ is an ideal triangulation of a compact manifold $M(T)$ then the Dehn filling of $D(T)$
along all the longitudes is the orientable double $D(M(T))$ of $M(T)$. Moreover
the cores of the filling solid tori are the doubles of the edges of $T$.}
\end{rem}

Proposition~\ref{homology:prop} now implies the following:

\begin{cor} \label{Haken:cor}
If in $T$ there are at most $n$ edges then
every Dehn filling of $D(T)$ has positive first Betti number. In
particular, it contains an incompressible surface.
\end{cor}

\begin{proof}
The assumption means that $D(T)$ has at most $n$ cusps.
A Dehn filling can decrease the first
Betti number by at most one, whence the result.
\end{proof}

Recall that the meridians of $D(T)$ were defined above. We now call \emph{non-meridinal}
a Dehn filling when it is non-meridinal on all cusps.
The following is a consequence of the results of Agol and Lackenby on Dehn fillings.

\begin{prop} \label{D:6:prop}
Let $V$ be a non-meridinal Dehn filling of $D(T)$.
\begin{enumerate}
\item If every edge of $T$ has valence at least $6$ then $V$ is Haken and the core of each filling
solid torus has infinite order in $\pi_1(V)$.
\item If every edge of $T$ has valence at least $7$ then $V$ is hyperbolic.
\end{enumerate}
\end{prop}

\begin{proof}
By Corollary~\ref{geometry:of:cusps:cor}, a non-meridinal slope on a cusp corresponding to
an edge with valence $q$ has length at least $q$ in the boundary of the horospherical cusp
neighbourhood $\calO(T)$.
Assume now all valences are at least $6$. Then the filled loops have length at least $6$ in $\calO(T)$,
so the proof of~\cite[Theorem 3.1]{lackenby}
implies that $V$ is irreducible and
$\partial$-irreducible, and the cores of the filling tori
have infinite order in
$\pi_1(V)$. Since a tetrahedron has $6$ edges, there can be at most $n$ edges in $T$, so
$V$ is Haken by Corollary~\ref{Haken:cor}. This proves the first assertion.
The second one now follows from the more general 6-theorem of Agol and Lackenby, together with Thurston's
hyperbolization theorem for Haken manifolds.
\end{proof}

\begin{rem} {\em
The set $\calD$ contains infinitely many examples,
including Agol's one~\cite{Agol:length-6},
which prove that the Agol-Lackenby theorem is sharp.
To see this, take $T$ with all edges having valence $6$, and assume
it is an ideal triangulation of an orientable manifold $M(T)$.
As noted in Remark~\ref{DM:rem}, the filling of $D(T)$ along the longitudes,
which have length $6$, is $D(M(T))$. But $D(M(T))$ is not hyperbolic, because
$M(T)$ is hyperbolic with cusps (a structure is obtained by
giving each tetrahedron of $T$ the shape of a regular ideal one),
so $D(M(T))$ contains essential tori.}
\end{rem}

\begin{prop} \label{T:6:prop}
Suppose that $T$ is an ideal triangulation of a manifold $M(T)$ and that all the edges of $T$
have valence at least 6. Then $M$ is hyperbolic
and no edge of $T$ is homotopic with fixed ends to an arc contained in
$\partial M$. In particular, the edges of $T$ are properly homotopic to geodesics.
\end{prop}

\begin{proof}
We assign angles to the edges of the tetrahedra of $T$,
choosing the angle $2\pi/v$ for an edge whose valence in $T$ is $v$.
The sum of angles around an edge of $T$ is then of course always $2\pi$,
and the assumption implies that the sum of the three angles at a vertex
of a tetrahedron is always at most $\pi$. This choice of angles
gives $T$ a structure weaker than the structure of \emph{angled triangulation}
defined by Lackenby~\cite{lackenby}, because he requires the sums of angles
at the vertices to be precisely $\pi$. However one can check that the arguments
in the proofs of \cite[Propositions 4.4 and 4.5]{lackenby} extend \emph{verbatim}
to our triangulation with angles, and they imply that $M(T)$ contains
no essential surface with non-negative Euler characteristic, so it is hyperbolic.
By Remark~\ref{DM:rem} and Proposition~\ref{D:6:prop}, the doubles
in $D(M(T))$ of the edges of $T$ are homotopically non-trivial, so
the edges of $T$ cannot be homotopic to arcs in the boundary.
\end{proof}

\paragraph{A closed incompressible surface}
An interesting property of $\cal D$ is described in the
following:

\begin{teo}\label{haken:teo}
Let $T$ be a triangulation which consists of at least two tetrahedra. Then $D(T)$
contains a closed incompressible surface of genus $2$ in its interior.
\end{teo}

\begin{proof}
Let $F$ be a face adjacent to two distinct tetrahedra $\Delta_1$
and $\Delta_2$ of $T$. Let $\Sigma_1$ and $\Sigma_2$ be two
(non-geodesic) parallel copies of the thrice-punctured sphere in
$D(T)$ corresponding to $F$, each $\Sigma_i$ contained inside the
block $K$ corresponding to $\Delta_i$. The surfaces $\Sigma_1$ and
$\Sigma_2$ cut the cusps they are incident to into $6$ annuli,
three small ones between them and three other ones $A_1$, $A_2$,
and $A_3$. We then define $\Sigma$ as $\Sigma_1\cup\Sigma_2\cup
A_1\cup A_2\cup A_3$. Since $\Sigma$ can be pushed in the interior of
$D(T)$, it is sufficient to show that it is incompressible.

Suppose by contradiction that $\Sigma$ is compressible in $D(T)$. Let us choose among all the
compressing discs for $\Sigma$ which are transverse to the
thrice-punctured spheres $S_i$ corresponding to the faces of $T$
the one, say $D$, whose intersection with $\bigcup S_i$
consists of the minimal number $k$ of components. We claim that $k=0$.
If $D\cap(\bigcup S_i)$ has a component which is a loop then the loop
must be trivial in one of the $S_i$'s (because they are geodesic surfaces). An easy
innermost argument then contradicts the minimality of $k$.

Suppose then that $D\cap(\bigcup S_i)$ has arc components, and choose an outermost one,
so there is an arc $\beta\subset\partial D$ such that $\alpha\cup\beta=\partial D'$, with
$D'$ a disc contained in $D$. Note that
$D'$ is embedded in a block $K$ as shown in Fig.~\ref{tubesandpants:fig},
with $\alpha$ lying on one $S_i$ and $\beta$ on $\Sigma$.
By construction, the surface $\Sigma$ is cut by the $S_i$'s
into some tubes and two pair of pants $P_1$ and $P_2$, contained in
distinct $K$'s. Each such tube or pair of pants is isotopic to a tube in a cusp
or a pair of pants inside a $S_i$, so it is incompressible.

The arc $\beta$ lies in a tube or in a pair of pants $P_i$.
If $\beta$ is not essential in it,
then $D'$ can be isotoped away, against the minimality of $k$.
This is always the case when $\beta$ lies on a tube, because $\alpha\subset S_i$ and $\alpha\cup\beta$ is a
loop, so the ends of $\beta$ must be on the same boundary component.
Therefore $\beta$ lies on a $P_i$, and again the
ends of $\beta$ must be on the same boundary component of the tube.
It easily follows that $\alpha\cup\beta$ is not null-homologous in
$K$, but $\alpha\cup\beta = \partial D'$, whence a contradiction.
This implies that indeed $k=0$.

We have shown that if $\Sigma$ is compressible then the intersection of $\Sigma$ with some block $K$
is compressible within the block, but in the course of our argument we have also noticed that this
cannot be the case. We then get a contradiction, and the proof is complete.\end{proof}

\section{Isometry groups}\label{isom:section}

Kojima proved in~\cite{Kojima:isom} that any finite group
is the full isometry group of a closed hyperbolic $3$-manifold.
We give here a new proof of this result, also providing an
upper bound on the volume of the manifold in terms
of the order of the group:

\begin{teo}\label{isom:teo}
There exists $k>0$ such that the following holds:
for any finite group $G$ there exists a closed orientable
hyperbolic $3$-manifold $V_G$ with
$\mathrm{Isom}(V_G)=\mathrm{Isom}^+(V_G)\cong G$
and $\mathrm{Vol}(V_G)\leqslant
k\cdot |G|^9$.
\end{teo}

Before plunging into the details of the proof we briefly describe the
scheme of our argument. In the whole section we fix a finite
group $G$ and we denote by $n$ its order. We first construct a
special polyhedron $Q_G$ with $\pi_1(Q_G)\cong G$. We carefully
choose a $Q_G$ having no symmetries, in such a way that the
group of combinatorial automorphisms of its universal cover
$P_G$ does not exceed the group of deck
transformations, which is isomorphic to $G$. Now let
$T_G$ be the triangulation dual to
$P_G$. Corollary~\ref{isometrie:cor} ensures that
$\mathrm{Isom}^+(D(T_G))\cong G$. In order to kill
the orientation-reversing isometries of $D(T_G)$
we then perform on $D(T_G)$ a suitable Dehn
filling, which  finally gives the desired $V_G$. Since volume
decreases under Dehn filling and the volume of $D(T_G)$ is proportional
to the number $c(P_G)$ of vertices of $P_G$, the bound on the
volume follows from a bound on $c(P_G)$.

\paragraph{Construction of the polyhedron}
We begin with the following:
\begin{lemma}\label{standpol:lem}
There exists $k>0$ independent of $G$ and a
special polyhedron $Q_G$ with $\pi_1(Q_G)\cong G$
and $c(Q_G)\leqslant k\cdot n^4$.
\end{lemma}

\begin{proof}
For the trivial group we can take a special spine of $S^3$, so we assume $n\geqslant 2$.
We consider the trivial presentation
for $G$ having as generators all the elements
of $G\setminus\{1\}$ and as relations all the expressions $abc^{-1}$ if
$a,b\neq 1$ and $a\cdot b=c\neq 1$ in $G$, and all the expressions $ab$ if $a\neq 1$ and $a\cdot b=1$ in $G$.
Let $U_n$ be the polyhedron obtained by
adding $n-1$ handles to $S^2$, thus getting a surface $\Sigma_n$,
and then attaching one disc to the core of each handle.
We now fix a bijective correspondence
between $G\setminus\{1\}$ and the handles, getting
an isomorphism between the free group
on $G\setminus\{1\}$ and $\pi_1(U_{n})$.

The relations in the presentation of $G$ translate
into simple loops $\gamma_1,\ldots,$ $\gamma_{(n-1)^2}$
on $\Sigma_{n}$. We
choose the $\gamma_i$'s in generic position with respect to each other and to the cores of the handles, and
intersecting each other and the cores in a minimal number $\nu$ of points.
The first condition and the fact that $n\geqslant 2$ easily imply that by attaching discs to $U_n$ along
the $\gamma_i$'s we get a special polyhedron $Q_G$ with $\pi_1(Q_G)\cong G$
and $\nu$ vertices.
The $\gamma_i$'s run $3n-4$ times (with multiplicity)
along each handle of $\Sigma_n$, and they give rise to $(n-1)(3n-4)$
properly embedded arcs in the central punctured sphere of $U_n$.
So $(n-1)(3n-4)$
vertices of $Q_G$ arise along the handles of $\Sigma_n$.
The other vertices lie in the central sphere of $U_n$, and there can be at most
$(n-1)(3n-4)(3n^2-7n+3)/2$ of them. Therefore
$\nu\leqslant k\cdot n^4$ for some $k>0$ independent of $G$.\end{proof}

\begin{prop}\label{asym:prop}
There exists $k>0$ independent of $G$ and
a special polyhedron $P_G$ with
group of combinatorial automorphisms isomorphic to $G$
and $c(P_G)\leqslant k\cdot n^9$.
\end{prop}

\begin{proof}
We suitably modify the polyhedron $Q_G$ given by the previous lemma
to ensure that its universal cover has the desired properties.

We begin with some definitions. If $P$ is a special
polyhedron, we say that
a vertex $v$ of $P$ is \emph{good} if
there is no edge of $P$ with both ends at $v$.
If $v$ is a good vertex of $P$ and $e$ is an
edge incident to $v$, we define a non-negative integer
$\lambda(v,e)$ as the maximal length of a simple simplicial path in the singular set of $P$
starting at $v$ with direction $e$ and touching,
apart from $v$, bad vertices only.

We first modify $Q_G$ by performing the local move described in
Fig.~\ref{bubble:fig} along each edge of $Q_G$.
\begin{figure}
\begin{center}
\mettifig{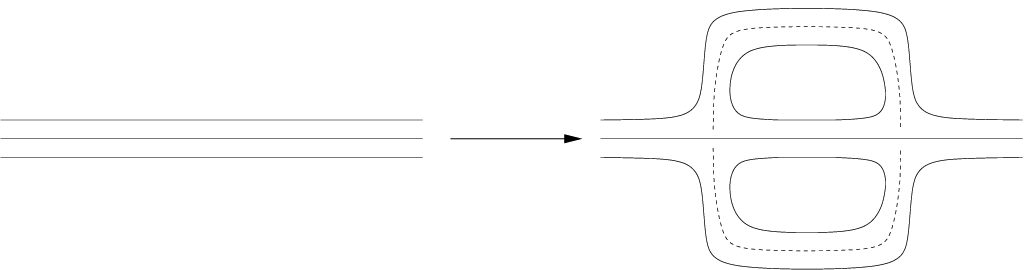}
\nota{Performing this move along each edge we eliminate bad vertices. The picture shows a neighbourhood
of the singular set, which determines a special polyhedron uniquely.} \label{bubble:fig}
\end{center}
\end{figure}
In this way we obtain a
polyhedron $Q'_G$ without bad vertices and such that
$\pi_1(Q'_G)\cong \pi_1(Q_G)$ and
$c(Q'_G)=5\cdot c(Q_G)$.
We now call \emph{curl} in a special polyhedron $P$
a face of $P$ whose
boundary consists of a single (closed) edge of $P$.
Let $(e_1,\ldots,e_\nu)$ be any ordering of the edges of $Q'_G$,
so $\nu=2\cdot c(Q'_G)$, and modify
$Q'_G$ by adding $i$ curls along $e_i$ as described in Fig.~\ref{curl:fig},
thus obtaining a polyhedron $Q''_G$.
\begin{figure}
\begin{center}
\mettifig{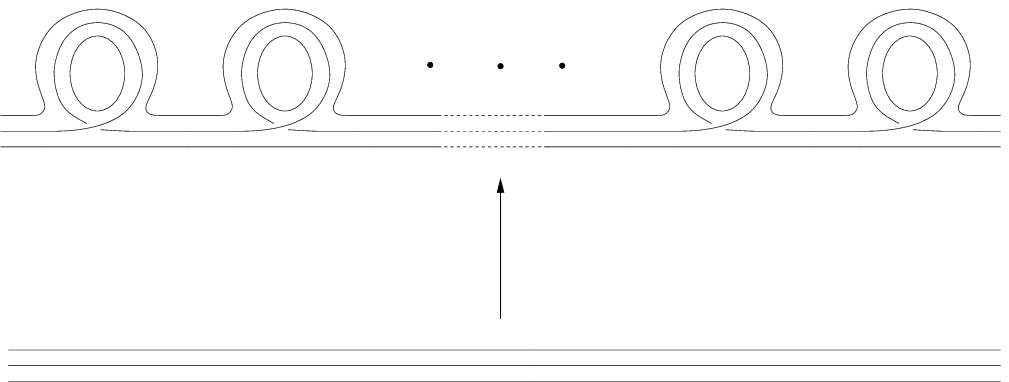}
\nota{Adding some curls along an edge of a special polyhedron.} \label{curl:fig}
\end{center}
\end{figure}
Observe that $Q''_G$
is homotopically
equivalent to $Q'_G$ and that $c(Q''_G)=2c(Q'_G)(c(Q'_G)+1)$.
We now define $P_G$ to be the universal cover of $Q''_G$.
The upper bound on $c(P_G)$ being obvious, we are left to prove the
first required property.

Let $\pi:P_G\to Q''_G$ be a fixed covering projection.
We first observe that the bad vertices of $Q''_G$
are precisely the $c(Q'_G)(2c(Q'_G)+1)$
vertices added to $Q'_G$. Moreover,
any edge $f$ starting from a good vertex $v$ of $Q''_G$ comes from a unique
edge $e_{i(f)}$ of $Q'_G$. By construction
we have $i(f)=\#\{\mathrm{curls\ along}\ e_{i(f)}\}=
\lambda(v,f)$.
The curls added to $Q'_G$ to get $Q''_G$ are contractible, so they
isomorphically lift
to $P_G$. This easily implies that a vertex of
$P_G$ is good if and only if its projection in $Q''_G$ is.
Moreover, if $\widetilde{v}\in\pi^{-1}(v)$
is a good vertex of $P_G$ and
$\widetilde{f}$ is
an edge starting from $\widetilde{v}$ with $\pi(\widetilde{f})=
f$, then $\lambda(\widetilde{v},\widetilde{f})=\lambda(v,f)=i(f)$.
Now let $\varphi$ be a combinatorial automorphism of $P_G$,
let $\widetilde{v}$ be a good vertex of $P_G$
and let $\widetilde{f}_1,\widetilde{f}_2,
\widetilde{f}_3,\widetilde{f}_4$ be the edges emanating
from $\widetilde{v}$.
Of course we have $\lambda(\varphi(\widetilde{v}),\varphi(\widetilde{f}_i))=
\lambda(\widetilde{v},\widetilde{f}_i)$ for $i=1,2,3,4$.
Since $Q'_G$ has no bad vertices and at least three good vertices,
this forces $\pi(\varphi(\widetilde{v}))=\pi(\widetilde{v})$
and $\pi(\varphi(\widetilde{f}_i))=\pi(\widetilde{f}_i)$
for $i=1,2,3,4$. It easily follows that
$\pi\compo\varphi=\pi$, \emph{i.e.} that $\varphi$ is a
deck transformation.
\end{proof}

\paragraph{Dehn filling}
Let $T_G$ be the triangulation dual to
the special polyhedron $P_G$ constructed during the proofs of
Lemma~\ref{standpol:lem} and Proposition~\ref{asym:prop}. The next result
concludes the proof of Theorem~\ref{isom:teo}.

\begin{prop}
There exists a hyperbolic Dehn filling $V_G$ of $D(T_G)$ with
$\mathrm{Isom}(V_G)=\mathrm{Isom}^+(V_G)\cong G$.
\end{prop}

\begin{proof}
We denote $N(T_G)$ by $N$ and
$D(T_G)$ by $D$.
Since $T_G$ has at least three good vertices,
it is not one of the exceptional
triangulations defined in Fig.~\ref{combinations:fig}. Then Corollary~\ref{isometrie:cor}
implies that $\mathrm{Isom}^+(D)\cong G$
and that there is a canonical
orientation-reversing involution $\tau$ of $D$.
Let $C_1,\ldots,C_m$ be the cusps of
$D$, and let $E_i$ be the boundary torus corresponding to
$C_i$ of the compactification of $D$. Note that $\tau$
extends to the $E_i$'s.
Since $\tau$ leaves invariant each $C_i$ and
reverses the orientation on it, there exists
exactly one slope $s_i$ on $E_i$ with $\tau(s_i)=s_i$.

Since $\mathrm{Aut}(P_G)$ is a group of covering
transformations and the 2-dimensio\-nal regions of $P_G$ are discs,
Brouwer's fixed point theorem implies that $\mathrm{Aut}(P_G)$ acts freely
on the regions, so
$\mathrm{Isom}^+(D)\cong
\mathrm{Isom}(N)\cong
\mathrm{Aut}(P_G)$
acts freely on the cusps. This implies that there exist
systems $(l_1,\ldots,l_m)$ of one slope per cusp which are invariant
under the action of $\mathrm{Isom}^+(D)$. More precisely,
Thurston's hyperbolic Dehn filling theorem implies that there
exists $(l_1,\ldots,l_m)$ with the following properties:
\begin{itemize}
\item[(A)]
$(l_1,\ldots,l_m)$ is $\mathrm{Isom}^+(D)$-invariant;
\item[(B)]
the filling $D(l_1,\ldots,l_m)$ of $D$ along the $l_i$'s is hyperbolic;
\item[(C)] the cores of the filling tori are the $m$ shortest geodesics of $D(l_1,\ldots,l_m)$;
\item[(D)] $l_i\neq s_i$ for $i=1,\ldots,m$.
\end{itemize}

If $h$ is an element of $\mathrm{Isom}^+(D)$,
from (A) we see that $h$ extends to an automorphism
of $D(l_1,\ldots,l_m)$, which in turn is homotopic
by Mostow's rigidity to an isometry $e(h)$.
It is easily seen that the map
$e:\mathrm{Isom}^+(D)\to
\mathrm{Isom}(D(l_1,\ldots,l_m))$ thus defined
is an injective group homomorphism.
To show it is surjective, pick $g\in\mathrm{Isom}(D(l_1,\ldots,l_m))$.
By (C), the cores of the filling tori are $g$-invariant, so $g$ restricts to an
automorphism of $D$, which by rigidity is homotopic to an isometry
$r(g)$. Moreover $r(g)$ leaves $(l_1,\ldots,l_m)$ invariant.
Now $r(g)=\tau^\varepsilon\circ h$
for some $h\in\mathrm{Isom}^+(D)$ and $\varepsilon\in\{0,1\}$.
If $r(g)$ maps $C_1$ to $C_i$ then $r(g)(l_1)=l_i$, but we know $h(l_1)=l_i$ and
$\tau(l_i)\neq l_i$ by (D), so $\varepsilon=0$.
It follows that $r(g)=h\in\mathrm{Isom}^+(D)$, so $g=e(h)$.
\end{proof}

\appendix
\section{Four-manifolds}
For the sake of completeness we recall here some facts proved 
in~\cite{Co:Th} concerning the class $\calD$,
even if they are not strictly speaking
needed for the present paper.

To begin we describe three
more spaces $\Theta(T,B)$
associated to a triangulation $T$ by the choice of a block $B$ 
as explained in Section~\ref{link:section}.  The third object is 
a 4-manifold $X(T)$ whose connection with $D(T)$ will be explained in the subsequent 
proposition.

\begin{itemize}
\item Take as a block the set $Q$ shown in
Fig.~\ref{blocchi:fig}-left, with the four \textsf{T}'s as faces.
Note that $Q$ naturally embeds in the polyhedron of Fig.~\ref{dualspine:fig}-right.
The result $\Theta(T,Q)$, that we denote
by $P_0(T)$, is a \emph{special polyhedron with boundary}.
\begin{figure}
\begin{center}
\mettifig{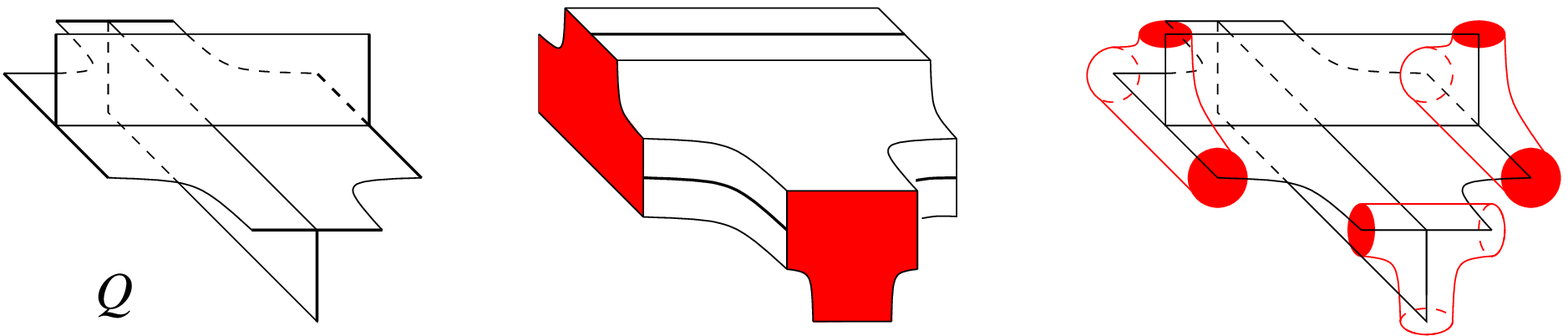,width=10 cm} \nota{The three blocks giving
respectively $P_0(T)$, $N(T)$, and $X(T)$.} \label{blocchi:fig}
\end{center}
\end{figure}
The boundary $\partial P_0(T)$
consists of loops that correspond to the edges of $T$.
We note that $P_0(T)$ can be identified to a regular neighbourhood
of the singular set of the special polyhedron $P(T)$ dual to $T$,
defined in Section~\ref{link:section}.
So $P(T)$ is obtained from $P_0(T)$ by attaching a disc to each
component of $\partial P_0(T)$;
\item Let $B$ be the 3-dimensional thickening
of $Q$ shown in Fig.~\ref{blocchi:fig}-centre, whose boundary contains four
shadowed hexagons (the faces) and $6$ arcs (three of which are visible in the picture).
Then $\Theta(T,B)$ is the relative handlebody $(H,\Gamma)=N(T)$, with loops $\Gamma$
constructed by attaching all the arcs.
Note that $P_0(T)\subset H$ and $\Gamma = \partial P_0(T)$, and that
$H$ collapses onto $P_0(T)$;
\item Take
as a block $B$ the $4$-dimensional thickening of $Q$,
given by the above $3$-dimensional thickening times the interval $[-1,1]$,
with faces defined as the products of the above hexagons with $[-1,1]$. 
Assign an arbitrary orientation to $B$ and note that
in the $3$-dimensional setting
each face-pairing induces an identification between hexagons, which may preserve or reverse
the orientation. To define the corresponding gluing of the faces
of the 4-dimensional blocks we then add to the gluing of the hexagons either the identity
or minus the identity of $[-1,1]$, so that the result is always
\emph{orientation-reversing}.
The space $\Theta(T,B)$, that we denote by $X(T)$,
is an \emph{oriented} $4$-dimensional manifold which contains $P_0(T)$
and collapses onto it.
\end{itemize}

\begin{prop} \label{X:prop}
\begin{itemize}
\item $\partial X(T)$ is a connected sum of copies of $S^2\times S^1$;
\item $\partial P_0(T)$ is a link in $\partial X(T)$;
\item $D(T)$ is the complement of $\partial P_0(T)$ in $\partial X(T)$.
\end{itemize}
\end{prop}

\begin{figure}
\begin{center}
\mettifig{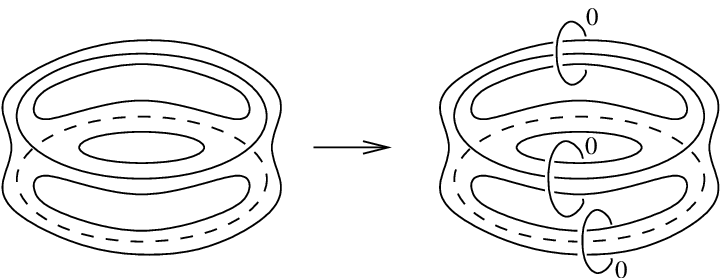,width=8cm}
\nota{From a description of $T$ via an immersion of $P_0(T)\subset S^3$ to a description of $D(T)$
via Dehn surgery.}
\label{link:fig}
\end{center}
\end{figure}

\begin{rem} {\em
The manifold $D(T)$ is not a link complement in $S^3$ in general. However,
as noted in~\cite{Co:Th}, we can always describe $D(T)$ as a link complement in $S^3$
with some $0$-surgered unknots. To do this, we
immerse $P_0(T)$ in $S^3$ as in the example of Fig.~\ref{link:fig}-left, we
take the link $\partial P_0(T)$,
and we encircle with $0$-surgered unknots
the triples of strands corresponding
the complement of a maximal tree in the singular set of $P_0(T)$.}
\end{rem}

The polyhedron $P_0(T)$ is a \emph{shadow} of the pair $(\partial
X(T),\partial P_0(T))$. Shadows are defined
in~\cite{Turaev:quantum} for arbitrary links in 3-manifolds, and
they can be used to show in
particular that $\calD$ is universal. 
Indeed the following holds (see~\cite{Turaev:quantum,Co:Th}):

\begin{teo}\label{universality:teo}
Every closed
orientable $3$-manifold is a Dehn
filling of $D(T)$ for some $T\in \calT$.
\end{teo}

This
theorem was used in~\cite{Co:Th} as a starting point to
define the \emph{shadow-complexity} of a $3$-manifold $Y$ as the
minimal number of vertices of a shadow for $Y$. The
shadow-complexity of $Y$ turns out to be strictly related to its
geometry.

\vspace{1.5 cm}

\noindent
Scuola Normale Superiore\\
Piazza dei Cavalieri 7, 56127 Pisa, Italy\\
f.costantino@sns.it\\
\vspace{.5 cm}

\noindent
Dipartimento di Matematica, Universit\`a di Pisa\\
Via F. Buonarroti 2, 56127 Pisa, Italy\\
frigerio@mail.dm.unipi.it, martelli@mail.dm.unipi.it\\
\vspace{.5 cm}

\noindent
Dipartimento di Matematica Applicata, Universit\`a di Pisa\\
Via Bonanno Pisano 25B, 56126 Pisa, Italy\\
petronio@dm.unipi.it

\end{document}